\input amstex
\magnification 1200
\vcorrection{-1cm}
\documentstyle{amsppt}
\NoBlackBoxes

\input epsf

\def\ODD  {{\lambda_{\operatorname{odd }}}}
\def\EVEN {{\lambda_{\operatorname{even}}}}
\def\odd {{\operatorname{odd}}}
\def\even {{\operatorname{even}}}

\def\Z {\bold Z }

\def\P {\bold P }

\def\R {\bold R }
\def\RP{\bold{RP}}

\def\calJ {{\Cal J}}
\def\calA {{\Cal A}}
\def\calF {{\Cal F}}
\def\union{\cup}

\def\<{\langle}
\def\>{\rangle}

\def\const{\operatorname{const}}
\def\Int{\operatorname{Int}}

\def\sign{\operatorname{sign}}

\def\Re{\operatorname{Re}}
\def\lk{\operatorname{lk}}
\def\card{\operatorname{card}}
\def\diag{\operatorname{diag}}
\def\Sign{\operatorname{Sign}}
\def\Null{\operatorname{Null}}
\def\eps{\varepsilon}

%=========================================================
\topmatter
%-----
\title
                Plane real algebraic curves of odd degree with a deep nest
\endtitle
%-----
\author
                S.Yu.~Orevkov
\endauthor
%-----
% \address        Steklov Mathematical Institute, 
%                 Russian Academy of Sciences,
%                 Gubkina 8, Moscow, 117966, Russia
% \endaddress
%-----
\email          orevkov\@mi.ras.ru
\endemail
%-----
\address        Laboratoire E.~Picard, UFR MIG,
                Universit\'e Paul Sabatier,
		118 route de Narbonne, 31062, Toulouse, France
\endaddress
%-----
\email          orevkov\@picard.ups-tlse.fr
\endemail
%-----
\abstract
		We apply the Murasugi-Tristram inequality to real algebraic
		curves of odd degree in $\RP^2$ with a deep nest, 
		i.e. a nest of the depth $k-1$ where $2k+1$ is the degree.
		For such curves, the ingredients of the Murasugi-Tristram
		inequality can be computed (or estimated) inductively
		using the computations for iterated torus links 
		due to Eisenbud and Neumann as the base case of the 
		induction and Conway's skein relation as the induction step.

		As an example of applications, we prove that some isotopy 
		types are not realizable by $M$-curves of degree 9.

		In Appendix B, we give some generalization 
		of the skein relation.
\endabstract
%-----
\endtopmatter
%==================================================================

\document

%===============================================
%==================================== References

\def\refCimasoni	{1} %{C}
\def\refEN              {2} %{EN}
\def\refFiedler		{3} %{F}
\def\refKauffman        {4} %{K}
\def\refNeumannSign     {5} %{N1}
\def\refNeumannJKTR     {6} %{N2}
\def\refOrevkovTop      {7} %{O1}
\def\refOrevkovFA       {8} %{O2}
\def\refOrevkovGAFA     {9} %{O3}

%===========================================
%================================== Sections

\def\sectDegreeNine	{2}
\def\sectSkein		{3}
\def\sectSplice		{4}
\def\sectMainLemma	{5}
\def\sectDet		{6}
\def\sectCycPol		{7}
\def\sectSign		{8}
\def\sectProofMainTh	{9}

%====================================================
%================================= Displaied formulas

\def\MainIneqOdd	{1}
\def\MainIneqEven	{2}
\def\defEpsnk		{3}
\def\defEpsnkPrime	{4}
\def\ineqOrient		{5}
\def\DegNineMainIneq	{6}
\def\eqDegNineCO	{7}
\def\eqDegNineNewCO	{8}
\def\ineqDegNineOrient	{9}
\def\ineqDegNineOriCor	{10}

\def\skein		{11}
\def\skeindet		{12}
\def\skeinsigma		{13}

\def\eqConway		{14}

\def\nonfib		{15}
				% bnkJ
\def\defdnkJ		{16}
\def\eqOmegaL		{17}

				% Appendix
\def\identityCycPol	{18}

\def\ineqMT		{19}
\def\eqNullZero		{20}
\def\ineqAodd		{21}
\def\eqAodd		{22}

				% Sect. Generalized skein relation
\def\genskein		{23}
\def\genskeindet	{24}
\def\genskeinDemoA	{25}
\def\genskeinDemoB	{26}
\def\genskeinDemoC	{27}
\def\genskeinD		{28}
\def\genskeindetD	{29}

%=======================================================
%=============================================== Figures

\def\figDegNineOne	{1}
\def\figDegNineTwo	{2}
\def\figCubic		{3}
\def\figSkein           {4}
\def\figCoreRemoved	{10}    
\def\figCoreRemained	{11}
\def\figlij		{12} 
\def\figTI              {13}
\def\figSD              {14}
\def\figLemmaA          {15}
\def\figLemmaB          {16}
\def\figBraid           {17}
\def\figBraidSD         {18}
\def\figBraidSDc        {19}
\def\figBaseSD		{20}

\def\figAodd		{21}
\def\figModif		{22}

\def\figSeifertSurface  {23}          
\def\figSeifertMatrix   {24}

%========================================================
%============================================== Proclaims

\def\MainTheorem	{1.1}

				% Application for deg 9
\def\thDegreeSeven	{2.1}
\def\lemDegNineBezout	{2.2}
\def\lemDegNineJump	{2.3}
\def\thDegNineOne	{2.4}
\def\corDegNine		{2.5}
\def\thDegNineTwo	{2.6}
\def\remDegNineOne	{2.7}
\def\remDegNineTwo	{2.8}

				% skein rel
\def\lemSkeinSigma	{3.1}
\def\corSkeinSigma	{3.2}

				% Splice diagrams

\def\thEN		{4.2}
\def\corEN		{4.3}
\def\remCimasoni	{4.4}
\def\lemTI              {4.5}

\def\MainLemma          {5.1}

\def\notBraid		{6.1}
\def\lemDet             {6.2}
\def\lemDetBase		{6.3}

\def\defCycPol		{7.1}
\def\lemDetCP           {7.2}
\def\lemCycPolA         {7.3}
\def\lemDetB		{7.4}
\def\lemDetA		{7.5}
\def\corCycPol		{7.6}
\def\lemSignDetA	{7.7}
\def\corBaseCycPol	{7.8}
\def\corDnkJ		{7.9}

\def\lemSignDelta	{8.1}
\def\propSignB		{8.2}
\def\propSignC		{8.3}

\def\lemCycPolB         {A.1}
\def\lemCycPolC         {A.2}
\def\lemCycPolD         {A.3}
\def\lemCycPolE         {A.4}

				% Gen skein rel

\def\propGenSkein       {B.2}
\def\corGenSkein        {B.3}

\head 1. Introduction and statement of the results
\endhead

In this paper we apply the Murasugi-Tristram inequality 
(as in [\refOrevkovTop--\refOrevkovGAFA])
to study real algebraic curves of odd degree in $\RP^2$ with a deep nest,
i.e. curves which have a nest of the depth $k-1$ where $2k+1$ is the 
degree. We study also analogous curves on real smooth ruled surfaces
(curves satisfying conditions (1)--(4) below). 
For curves with a deep nest, the braid defined in [\refOrevkovTop] is uniquely 
determined by the arrangement of the curve with respect to the pencil
of lines centered at a point inside the nest
(the arrangement with respect to the fibers on the ruled surfaces).
Moreover, if the 
degree is odd then the right hand side of the Murasugi-Tristram inequality
(the signature and the nullity of the braid) can be computed 
inductively in some cases (and estimated in the other cases).
The base case of the induction 
uses the computations for iterated torus links due to Eisenbud and Neumann,
and the induction step is Conway's skein relation.
As an example of applications, we prove that some isotopy types are
not realizable by $M$-curves of degree 9.

In Appendix B, we give some generalization of the skein relation.

Let us fix an integer $n\ge 1$.
Let $\pi_n:\Sigma_n\to\P^1$ be the
fiberwise compactification of the line bundle ${\Cal O}(n)$. 
Let $E_n$ be the infinite section (i.e. $E_n^2=-n$). 
If $n=1$ then $\Sigma_n$ is 
the blown up $\P^2$. Otherwise $\Sigma_n$ is a minimal
rational smooth ruled surface --- Hirzebruch surface.
If $n\ge 1$  then $\Sigma_n$ has a unique real form 
($\Sigma_0$ has two real forms: hyperboloid and ellipsoid
but we shall not consider the case $n=0$). 
For a real algebraic variety $X$, 
we denote the set of its real (resp. complex) points by $\R X$
(resp. by $X$).

We have a fibration
$\pi_n:\R \Sigma_n\to\RP^1$ with the fiber $\RP^1$.
If $n$ is even then $\R\Sigma_n$ is a torus; 
if $n$ is odd then $\R\Sigma_n$ is a Klein bottle.

For an algebraic curve $A$ in $\Sigma_n$, we define the {\it bidegree of} $A$
as $(A\cdot F, A\cdot E_n)$ where $F$ is a fiber of $\pi_n$.
The curves of a given bidegree $(m,l)$ form a linear system on $\Sigma_n$.
If $(x,y)$ is a coordinate system in $\Sigma_n\setminus(F\cup E_n)$ such that
the fibers are given by $\{x=\const\}$ then the Newton polygon
of a generic curve of bidegree $(m,l)$ is
$(0,0)$-$(l+mn,0)$-$(l,m)$-$(0,m)$.
The genus of such a curve is $g_{m,l}=(m-1)(mn/2+l-1)$ 
(the number of interior integral points in the Newton polygon).  

In this paper we study real algebraic curves $A$ on 
$\Sigma_n$ satisfying the following conditions:

\roster
\item
    $A$ is non-degenerate;
\item
    the bidegree of $A$ is $(m,0)$ where $m=2k+1$ is odd;
\item
    There exists $A_0\subset\R A$ such that $\pi_n|_{A_0}$ is a
    covering of $\RP^1$ of degree $m-2$;
\item
    The curve $A$ is not hyperbolic, i.e.
    for any connected component $V$ of $\R A\setminus A_0$, the 
    topological degree of the mapping $\pi_n|_V:V\to\RP^1$ is zero.
\endroster

\example{Example} If $n=1$ then $\Sigma_1$ is $\P^2$ 
blown-up at a real point $p$ and
$A$ is the strict transform of a smooth real curve of degree
$m=2k+1$ with a nest of ovals $O_1,\dots,O_{k-1}$, such that
$p\in\Int O_1\subset\dots\subset\Int O_{k-1}$.
\endexample

Suppose that a curve $A$ on $\Sigma_n$ satisfies the conditions (1)--(4).
If $n$ is odd then $A_0$ is a disjoint union of $k$ circles
$A_0=\calJ\sqcup O_1\sqcup O_2\sqcup\dots\sqcup O_{k-1}$, moreover,
$\pi_n|_{\calJ}:\calJ\to\RP^1$ is a diffeomorphism and
$\pi_n|_{O_i}:O_i\to\RP^1$ is a double covering.
We call $\calJ$ {\it the odd branch of $A$}.
%we call $O_1,\dots,O_{k-1}$ {\it the double branches of $A$}.
If $n$ is even then $A_0$ is a disjoint union of $2k$ circles,
each being mapped diffeomorphically onto $\RP^1$ by the projection
$\pi_n$. In this case we fix any component $\calJ$ of $A_0$ and call it 
the {\it odd branch}.
 
For any parity of $n$ we
call the connected components of $\R A\setminus A_0$ 
{\it the ovals of $A$}.
An oval $V$ of $A$ is called {\it odd} (resp. {\it even}) if a generic path 
in $\R\Sigma_n\setminus E_n$ relating $V$ to $\calJ$ intersects
$A_0$ in an odd (resp. even) number of points.

% Let us fix an orientation on $\RP^1$.
% Any to points $p,q\in\RP^1$ divide $\RP^1$ into two segments.
% We shall denote them by $[p,q]$ and $[q,p]$ according to the chosen
% orientaiton of $\RP^2$.

\definition{ Definition } Let $A$ be a nonsingular 
real algebraic curve on 
$\Sigma_n$ satisfying the conditions (1)--(4). 
A {\it jump over $\calJ$} is a pair of ovals $(V_1,V_2)$
such that for any $p_1\in V_1$, $p_2\in V_2$,
the ribbon 
% $D=\pi_n^{-1}\big([\pi_n(p_1),\pi_n(p_2)]\big)$ 
$D=\pi_n^{-1}([q_1,q_2])$, $q_j=\pi_n(p_j)$,
does not contain other ovals of $A$ and the points $p_1$ and $p_2$
belong to different connected components of 
$D\setminus(E_n\union\calJ)$
(here we suppose that an orientation is fixed on $\RP^1$ and 
the segment $[q_1,q_2]$ is oriented from $q_1$ to $q_2$).
\enddefinition

It is clear that the number of jumps is of the same parity as $n$.

A real algebraic curve $A$ is called an {\it $M$-curve} if $\R A$ has 
$g+1$ connected components where $g$ is the genus of $A$.
It is called an $(M-r)$-curve if $\R A$ has $g+1-r$ 
connected components.
If $A$ is a curve of bidegree $(m,0)$ on $\Sigma_n$ then 
$g=(m-1)(mn-2)/2$.

The main result of the paper is the following.

\proclaim{ Theorem \MainTheorem }
Let $A$ be a real algebraic $(M-r)$-curve of bidegree $(m,0)$, 
$m=2k+1$, on $\Sigma_n$ satisfying the conditions (1)--(4). 
If $n$ is divisible by $4$, we suppose that $k=1$.
Let $\calJ$ be chosen as above
and $J$ be the number of jumps over $\calJ$.
Let us denote 
the number of all ovals by $\lambda$, 
the number of odd ovals by $\ODD$, and
the number of even ovals by $\EVEN$.
We shall suppose that $\lambda>J>0$.
%
% and if $J\equiv(2k-1)n\mod4$
% then we suppose also that $\lambda>J$.

Then one has
$$
	|nk^2-3k+1+\eps_{k,n}-r-J| \le r+2\ODD+
	\cases 
		k-1, &\text{if $n$ is odd}\\
		2(k-1), &\text{if $n$ is even}
	\endcases 
							\eqno(\MainIneqOdd)
$$
%and
$$
	|nk^2-3k+1+\eps'_{k,n}-r-J| \le r+2\EVEN+
	\cases 
		k-1, &\text{if $n$ is odd}\\
		2(k-1), &\text{if $n$ is even}
	\endcases 
							\eqno(\MainIneqEven)
$$
where
$$
\eps_{n,k}={1+(-1)^k\over2}\cdot\Re i^{\,n-1} =
\cases
	(-1)^{(n-1)/2}	&\text{if $k+1\equiv n\equiv 1\mod2$,}\\
	0		&\text{otherwise}
\endcases
							\eqno(\defEpsnk)
$$
%and
$$
\eps'_{n,k}={1-3(-1)^k\over2}\cdot\Re i^{\,n-1} =
\cases
	(-1)^{(n+1)/2}		&\text{if $k+1\equiv n\equiv 1\mod2$,}\\
				%&\text{if $n$ is odd and $k$ is even,}\\
				%&\text{if $(n,k)\equiv(1,0)\mod2$}\\
	2\cdot(-1)^{(n-1)/2}	&\text{if $k\equiv n\equiv 1\mod2$,}\\
	0			&\text{otherwise}.
\endcases
							\eqno(\defEpsnkPrime)
$$
\endproclaim

This theorem will be proved in Section \sectProofMainTh.

Now let us suppose that the curve $A$ is {\it dividing}, 
or {\it of the type I}. This means that $A\setminus\R A$
consists of two halves which are mapped onto each other by the complex
conjugation. An orientation of $\R A$ induced by the complex orientaion
of one of the halves is called a {\it complex orientation}. There are 
two opposite complex orientation. Let us fix one of them.
An oval $V$ of $A$ is called {\it positive} (resp. {\it negative})
if $[v] = -2[\calJ]$ (resp. $[v] = 2[\calJ]$) in 
$H_1(\R\Sigma_n\setminus E_n)$.
Let us denote the number of positive (resp. negative) ovals by $\lambda_+$
(resp. by $\lambda_-$).
By Fiedler's alternating orientation rule [\refFiedler]
$$
	J \ge |\lambda_+ - \lambda_-|.				\eqno(\ineqOrient)
$$
Indeed, each sequence of ovals between successive jumps
over $\calJ$ contrubutes at most $1$ to the right hand side of
(\ineqOrient).

Combination of (\MainIneqOdd), (\ineqOrient), and the
formula for complex orientations from [\refOrevkovTop; Theorem 1.4A]
allows to exclude some isotopy types of $M$-curves of degree $9$
on $\RP^2$ (see Proposition \thDegNineOne).

%%%%%%%%%%%%%%%%%%%%%%%%%%%%%%%%%%%%%%%%%%%%%%%%%%%%%%%%%%%%%%%%%%%%%%
%%%%%%%%%%%%%%%%%%%%%%%%%%%%%%%%%%%%%%%%%%%%%%%%%%%%%%%%%%%%%%%%%%%%%%
%%%%%%%%%%%%%%%%%%%%%%%%%%%%%%%%%%%%%%%%%%%%%%%%%%%%%%%%%%%%%%%%%%%%%%
%%%%%%%%%%%%%%%%%%%%%%%%%%%%%%%%%%%%%%%%%%%%%%%%%%%%%%%%%%%%%%%%%%%%%%
%%%%%%%%%%%%%%%%%%%%%%%%%%%%%%%%%%%%%%%%%%%%%%%%%%%%%%%%%%%%%%%%%%%%%%
%%%%%%%%%%%%%%%%%%%%%%%%%%%%%%%%%%%%%%%%%%%%%%%%%%%%%%%%%%%%%%%%%%%%%%
%%%%%%%%%%%%%%%%%%%%%%%%%%%%%%%%%%%%%%%%%%%%%%%%%%%%%%%%%%%%%%%%%%%%%%

\head \sectDegreeNine. 
		Application for $M$-curves of degrees $7$ and $9$ on $\RP^2$
\endhead
In this section we give details of the proofs of the results of 
[\refOrevkovFA].

\subhead \sectDegreeNine.1. Curves of degree $7$
\endsubhead

\proclaim{ Theorem \thDegreeSeven }{\rm[\refOrevkovFA; Theorem 3].}
There does not exist an $M$-curve of degree $7$ on $\RP^2$
with the complex scheme 
$\<\calJ \sqcup 10_+ \sqcup 3_- \sqcup 1_-\<1_-\>\>$.
\endproclaim

\demo{ Proof }
Blowing up a point inside the innermost oval, we are in the hypothesis
of Theorem \MainTheorem\ with $n=1$, $k=3$, $r=0$, $\eps_{n,k}=0$.
Thus, (\MainIneqOdd) yeilds $|J-1|<2$ which contradicts
to (\ineqOrient) because we have $\lambda_+=10$ and $\lambda_-=3$.
\qed
\enddemo

%%%%%%%%%%%%%%%%%%%%%%%%%%%%%%%%%%%%%%%%%%%%%%%%%%%%%%%%%%%%%%%%%%%%%%

\subhead \sectDegreeNine.2. Curves of degree $9$
\endsubhead
Let $A$ be an $M$-curve of degree $9$ on $\RP^2$ whose isotopy type
is $\<\calJ\sqcup \alpha \sqcup 1\< \beta \sqcup 1\< \gamma \>\>\>$
with $\gamma\ge1$
(see Figure \figDegNineOne).
Abusing notation, we shall denote the set of empty ovals in the 
corresponding region by $\<\alpha\>$, $\<\beta\>$, or $\<\gamma\>$.
Let us denote the odd branch by $\calJ$ and the 
non-empty ovals by $O_1$ and $O_2$.

\midinsert 
\epsfxsize 105mm
\centerline{\epsfbox{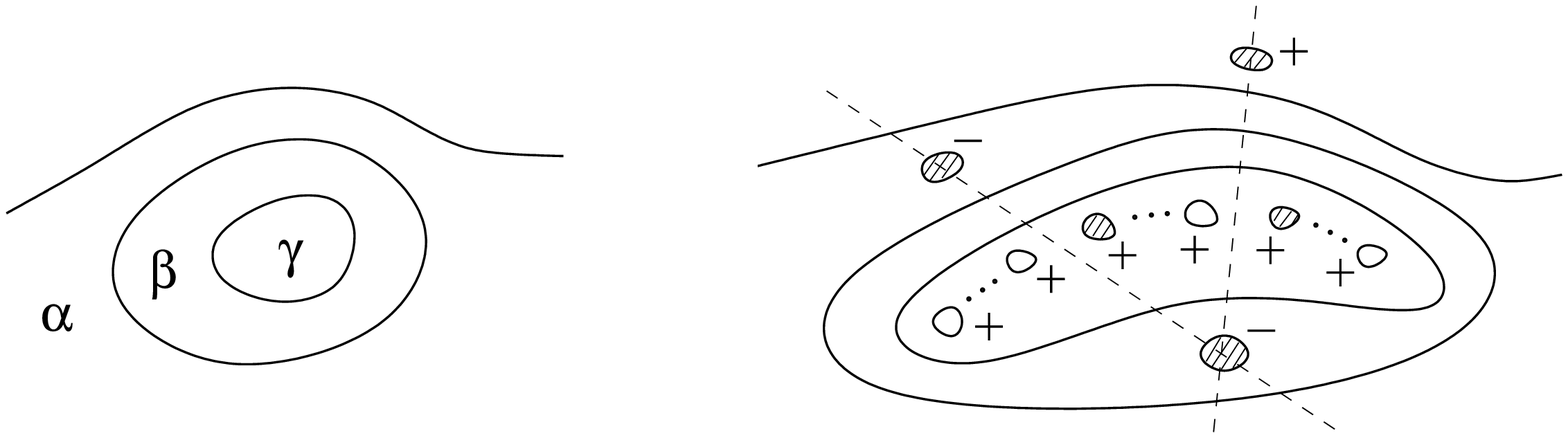}}
\botcaption{
		Fig.~\figDegNineOne \hbox to 50mm{} 
		Fig.~\figDegNineTwo \hbox to 10mm{} 
          } 
\endcaption
\endinsert

For $v\in\<\gamma\>$, let us denote by $J_v$
the number of jumps through $\calJ$ in the pencil of
lines centered at a point $p\in\Int v$
(i.e. the number of jumps for the transform of $A$ on $\RP^2$
blown up at $p$). Then we have $n=1$, $k=4$, $r=0$, $\eps_{n,k}=1$ and 
hence, (\MainIneqOdd) reads as
$$
	3-2\beta\le J_v \le 9 + 2\beta.
							\eqno(\DegNineMainIneq)
$$

\proclaim{ Lemma \lemDegNineBezout } 
If ovals $v_1,v_2\in\<\gamma\>$ are separated by a line 
through ovals $v_3,v_4\in\<\alpha\>$ then $J_v=1$ for any $v\in\<\gamma\>$.
\endproclaim

\demo{ Proof } Any conic through $v_1,\dots,v_4$ meets $O_1\cup O_2$
at $\ge 8$ points. Hence, if it meets one more empty oval, it cannot meet 
$\calJ$. Therefore, if we choose a point inside each empty oval then
these points are vertices of a convex polygon which does not intersect $\calJ$
and the result follows.
\qed
\enddemo

Let the complex scheme of $A$ be 
$\<\calJ\sqcup \alpha_+\sqcup\alpha_-\sqcup 
1_{\eps_2}\< \beta_+\sqcup\beta_-\sqcup 1_{\eps_1}\< \gamma \>\>\>$ 
where
$\eps_1,\eps_2\in\{\pm\}$. Let us set
$\delta\alpha=\alpha^+ - \alpha^-$, 
$\delta\beta =\beta^+  - \beta^- $, 
$\delta\gamma=\gamma^+ - \gamma^-$, 
The Rohlin-Mishachev's formula for complex orientations implies
$$
	\delta\alpha + \eps_2 + (1-2\eps_2)(\delta\beta + \eps_1) +
		(1-2\eps_1-2\eps_2)\delta\gamma = 8.
							\eqno(\eqDegNineCO)
$$
The formulas for complex orientations [\refOrevkovTop; Theorem 1.5A] imply
% $$
% 	{\eps_2+1\over2}(\delta\beta + \delta\gamma) + 
% 	{\eps_1+1\over2}\delta\gamma 
% 		= \left({\eps_2+1\over2}+{\eps_1+1\over2}\right)^2.
% 							\eqno(\eqDegNineNewCO)
% $$
$$
	(\eps_2+1)(\delta\beta + \delta\gamma) + 
	(\eps_1+1)\delta\gamma 
		= (\eps_2+\eps_1+2)^2/2.
							\eqno(\eqDegNineNewCO)
$$
The inequality (\ineqOrient) takes the form
$$
	J_v \ge |\delta\alpha + \delta\beta + \delta\gamma - \sign v|.
							\eqno(\ineqDegNineOrient)
$$
Combining it with (\DegNineMainIneq), we obtain
$$
	|\delta\alpha+\delta\beta+\delta\gamma\mp1|\le9+2\beta
	\qquad\text{if $\gamma_\pm\ge0$.}
							\eqno(\ineqDegNineOriCor)
$$

\proclaim{ Lemma \lemDegNineJump } 
If $|\delta\gamma|>1$ and $\alpha>0$ then $\beta>0$.
\endproclaim

\demo{ Proof } Suppose that $|\gamma_+ - \gamma_-|>1$, $\alpha>0$, and $\beta=0$.
Let us consider a pencil of lines through an oval from $\<\alpha\>$.
By Fiedler's alternating orientation rule, some two ovals from $\<\gamma\>$ 
must be separated in this pencil by a line through another oval from $\<\alpha\>$.
By Lemma \lemDegNineBezout, this implies $J_v=1$ which contradicts
(\DegNineMainIneq).
\qed
\enddemo

\proclaim{ Theorem \thDegNineOne } 
If an $M$-curve of degree nine is of the form
$\<\calJ\sqcup \alpha \sqcup 1\< 1\< \gamma \>\>\>$
where $\gamma>0$,
then its complex scheme is one of the following

\smallskip
\centerline{
$
\<\calJ\sqcup 
({\alpha+7\over 2})_+ \sqcup
({\alpha-7\over 2})_-
 \sqcup 1_-\< 1_-\< 
({\gamma+1\over 2})_+ \sqcup
({\gamma-1\over 2})_-\>\>\>,
$}

\smallskip
\centerline{
$
\<\calJ\sqcup 
({\alpha+7\over 2})_+ \sqcup
({\alpha-7\over 2})_-
 \sqcup 1_\pm\< 1_\mp\< 
({\gamma-1\over 2})_+ \sqcup
({\gamma+1\over 2})_-\>\>\>.
$}
\endproclaim

\demo{ Proof }
Only these complex schemes satisfy
%These are the only complex schemes satisfying 
(\eqDegNineCO), (\eqDegNineNewCO), (\ineqDegNineOriCor) and Lemma \lemDegNineJump.
\qed
\enddemo

\proclaim{ Corollary \corDegNine } {\rm[\refOrevkovFA; Theorems 1 and 2].}
If an $M$-curve of degree nine is of the form
$\<\calJ\sqcup \alpha \sqcup 1\< 1\< \gamma \>\>\>$, 
where $\gamma>0$,
then $\alpha$ is odd and $\alpha\ge 7$.
\qed\endproclaim

\proclaim{ Theorem \thDegNineTwo } The isotopy types
$\<\calJ \sqcup 2 \sqcup 1\<1\sqcup 1\<23\>\>\>$ and
$\<\calJ \sqcup 3 \sqcup 1\<1\sqcup 1\<22\>\>\>$
are not realizable by $M$-curves of degree nine.
\endproclaim

\demo{ Proof } 
% {\bf 1.} 
$\<\calJ \sqcup 2 \sqcup 1\<1\sqcup 1\<23\>\>\>$.
%(see [\refOrevkovTop; {\it Erratum}\,]).
It follows from (\eqDegNineCO) and  (\eqDegNineNewCO), the complex scheme must be
$\<\calJ \sqcup 1_+ \sqcup 1_- \sqcup 1_-\<1_-\sqcup 1_-\<13_+\sqcup 10_-\>\>\>$.
Applying Fiedler's orientation alternating rule for a pencil of lines
centered in the oval from $\<\beta\>$, we see that the curve must be arranged
with respect to some two lines as in Figure {\figDegNineTwo}.
Then a conic passing through the five shadowed ovals contradicts
Bezout's theorem.
\medskip

% {\bf 2.} 
$\<\calJ \sqcup 3 \sqcup 1\<1\sqcup 1\<22\>\>\>$.
It follows from (\eqDegNineCO) and  (\eqDegNineNewCO), 
that the complex scheme must be
$\<\calJ \sqcup 1_+ \sqcup 2_- \sqcup 1_-\<1_+\sqcup 1_-\<12_+\sqcup 10_-\>\>\>$.
Let us consider a pencil of lines centered in an oval from $\<\alpha\>$.
By Fiedler's orientation alternating rule, the ovals $\<\gamma\>$ must be 
separated by a line through an oval from $\<\beta\>\cup\<\alpha\>$. 
It cannot be the oval from $\<\beta\>$ because the latter is positive. 
Hence, the hypothesis of
Lemma \lemDegNineJump\ are satisfied, 
and we have $J_v=1$ for any $v\in\<\gamma\>$.
This contradicts to (\ineqDegNineOrient).
\qed
\enddemo

\remark{ Remark \remDegNineOne }
In [\refOrevkovTop; Corollary 1.8], we announced
the nonrealizability of $9$ isotopy types by $M$-curves of $9$-th degree.
The proof given in [\refOrevkovTop] was based on a wrong assertion of
[\refOrevkovTop; Lemma 1.9] (see [\refOrevkovTop; {\it Erratum}\,]).
However, Corollary {\corDegNine} and Theorem {\thDegNineTwo} imply
that the statement of [\refOrevkovTop; Corollary 1.8] was correct.
\endremark

\remark{ Remark \remDegNineTwo }
We give Theorem {\thDegNineTwo} not as an application of Theorem {\MainTheorem}
(which is not used in the proof) but just by the reasons explained 
in Remark \remDegNineOne.
\endremark

%%%%%%%%%%%%%%%%%%%%%%%%%%%%%%%%%%%%%%%%%%%%%%%%%%%%%%%%%%%%%%%%%%%%%%
%%%%%%%%%%%%%%%%%%%%%%%%%%%%%%%%%%%%%%%%%%%%%%%%%%%%%%%%%%%%%%%%%%%%%%
%%%%%%%%%%%%%%%%%%%%%%%%%%%%%%%%%%%%%%%%%%%%%%%%%%%%%%%%%%%%%%%%%%%%%%
%%%%%%%%%%%%%%%%%%%%%%%%%%%%%%%%%%%%%%%%%%%%%%%%%%%%%%%%%%%%%%%%%%%%%%
%%%%%%%%%%%%%%%%%%%%%%%%%%%%%%%%%%%%%%%%%%%%%%%%%%%%%%%%%%%%%%%%%%%%%%
%%%%%%%%%%%%%%%%%%%%%%%%%%%%%%%%%%%%%%%%%%%%%%%%%%%%%%%%%%%%%%%%%%%%%%
%%%%%%%%%%%%%%%%%%%%%%%%%%%%%%%%%%%%%%%%%%%%%%%%%%%%%%%%%%%%%%%%%%%%%%

\head \sectSkein. Signature, determinant and skein-relation
\endhead
Recall some definitions. Let $B$ be a real symmetric matrix and
$D=QBQ^T$ its diagonalization. {\it The signature} $\Sign B$
and {\it nullity} $\Null B$ of $B$ are defined as the sum of signs
of the diagonal entries and the number of zeros on the diagonal of $D$.

{\it A Seifert surface} of an oriented link $L$ in the 3-sphere
$S^3$ is a conncted oriented surface $X\subset S^3$ whose boundary
(taking into account the orientation) is $L$. Let $x_1,\dots,x_n$
be a base in $H_1(X,\Z)$.
{\it A Seifert matrix} is the matrix of the linking numbers
$\lk(x_i,x_j^+)$ where $x_j^+$ is the result of
a small shift of the cycle $x_j$ along a positive normal field to $X$.
The {\it signature} $\Sign L$ and the {\it nullity}
$\Null L$ of $L$ are defined by
\footnote{The nullity of $L$ is defined usually as $1+\Null(V+V^T)$.}
$\Sign L = \Sign(V+V^T)$, $\Null L = \Null(V+V^T)$ 
where $V$ is a Seifert matrix of $L$.

The {\it Conway potential function} of a link $L$ is defined 
as $\Omega_L(t)=\det( t^{-1}V - t V^T)$ 
where $V$ is a Seifert matrix of $L$.
By definition, $\Omega_L(t)=1$ if $L$ is the trivial knot.
We define the {\it determinant} of $L$ as 
$\det L=\Omega_L(i)$ where $i=\sqrt{-1}$.
It is known that $\Omega_L(t)$ is a link invariant and the following
{\it skein-relation} holds (see [\refKauffman]).
Let $L_-$, $L_0$, and $L_+$ be links whose diagrams coincide everywhere
except of some disk where they look as in Fig.~\figSkein. Then 
$$
    \Omega_{L_+}(t) - \Omega_{L_-}(t) = (t-t^{-1})\Omega_{L_0}(t).
                                                          \eqno(\skein)
$$
Substituting $t=i$ into (\skein), we obtain
$$
	\det L_+ - \det L_- = 2i\det L_0.		\eqno(\skeindet)
$$

\midinsert 
\epsfxsize 90mm
\centerline{\epsfbox{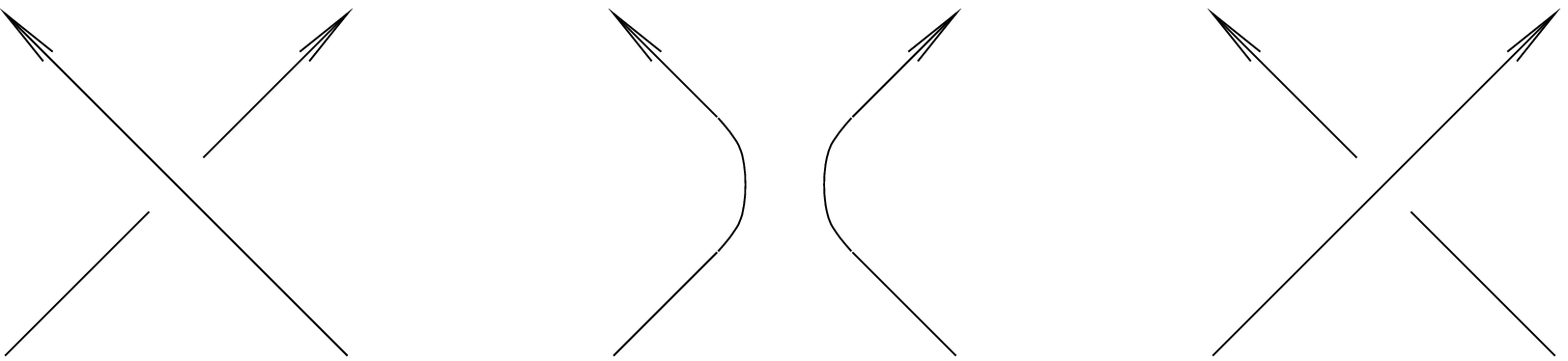}}
\centerline{$L_-$ \hbox to 28mm{} $L_0$ \hbox to 28mm{} $L_+$}
\botcaption{
             Fig.~\figSkein
          } 
\endcaption
\endinsert

Let $L$ and $L'$ be oriented links.
Say that $L'$ is obtained from $L$ by a {\it band-attachement}
if $L'=(L\setminus b(\partial I\times I))\cup b(I\times\partial I)$
where $I=[0,1]$ and $b:I^2\to S^3$ is an embedding such that 
$L\cup b(I^2) = b(\partial I\times I)$ and the orientations of
$L$ and $L'$ coincide on $L\cap L'$.
For instance, $L_+$ and $L_-$ in Fig.~\figSkein\ are obtained from $L_0$
by a band-attachement.
It is clear that if $L'$ is obtained from $L$ by a band-attachement then
$L$ also can be obtained from $L'$ by a band-attachement.

\proclaim{ Lemma \lemSkeinSigma } {\rm(Conway).}
Let a link $L'$ be obtained by a band-attachement 
from a link $L$. 

(a). If $\det L\ne 0$ then
$$
   \Sign L' = \Sign L + 
		%s
		%\quad\text{and}\quad
	   	%\Null L' = 1 - |s|
	   	%\qquad\text{where}\;\;
	   	%s = 
	\sign{i\det L' \over \det L}
                                                        \eqno(\skeinsigma)
$$

(b). $|\Null L' - \Null L| + |\Sign L' - \Sign L| = 1$. 
% Moreover, $|\Null L' - \Null L| = 1$ if and only if $\Sign L' = \Sign L$.
\endproclaim

\remark{ Remark } $\det L\in\Z$ (resp. $\in i\Z$)
if the number of components of $L$ is odd (resp. even)
and the band attachement changes the parity of the number of components.
Hence, the fraction in (\skeinsigma) is real.
\endremark

\demo{ Proof } Let $b:I\times I\to S^3$ be the attached band.
One can choose a Seifert surface $X$ for $L$
such that $X\cap b(I^2) = b(\partial I\times I)$. Then
$X'=X\cup b(I\times I)$ is a Seifert surface for $L'$
and one can choose a base of $H_1(X')$ by adding one element to
a base of $H_1(X)$.
Let $V$ and $V'$ be the corresponding Seifert matrices.
Then $V$ is a principal $n\times n$-minor of an 
$(n+1)\times(n+1)$-matrix $V'$. This implies the part (b).
To prove the part (a), note that
$$
 \Sign L' - \Sign L = \sign{\det(V'+{V'}^T)\over \det(V+V^T)}
  = \sign{ i^{n+1}\det L' \over i^{n}\det L}
$$
\qed\enddemo

\proclaim{ Corollary \corSkeinSigma } 
Let each of two links $L'$ and $L''$ be
obtained by a band-attachement from the same link $L$. 
If $\det L'=\det L''\ne 0$ then $\sigma(L'')=\sigma(L')$.\qed
\endproclaim

%%%%%%%%%%%%%%%%%%%%%%%%%%%%%%%%%%%%%%%%%%%%%%%%%%%%%%%%%%%%%%%%%%%%%%
%%%%%%%%%%%%%%%%%%%%%%%%%%%%%%%%%%%%%%%%%%%%%%%%%%%%%%%%%%%%%%%%%%%%%%
%%%%%%%%%%%%%%%%%%%%%%%%%%%%%%%%%%%%%%%%%%%%%%%%%%%%%%%%%%%%%%%%%%%%%%
%%%%%%%%%%%%%%%%%%%%%%%%%%%%%%%%%%%%%%%%%%%%%%%%%%%%%%%%%%%%%%%%%%%%%%
%%%%%%%%%%%%%%%%%%%%%%%%%%%%%%%%%%%%%%%%%%%%%%%%%%%%%%%%%%%%%%%%%%%%%%
%%%%%%%%%%%%%%%%%%%%%%%%%%%%%%%%%%%%%%%%%%%%%%%%%%%%%%%%%%%%%%%%%%%%%%
%%%%%%%%%%%%%%%%%%%%%%%%%%%%%%%%%%%%%%%%%%%%%%%%%%%%%%%%%%%%%%%%%%%%%%
%%%%%%%%%%%%%%%%%%%%%%%%%%%%%%%%%%%%%%%%%%%%%%%%%%%%%%%%%%%%%%%%%%%%%%

\head \sectSplice. Iterated torus links and their splice diagrams
\endhead
Let $L$ be an oriented link in $S^3$ and $S$ a component of $L$.
Let $T$ be a tubular neighbourhood of $S$ which 
does not intersect $L\setminus S$. Let $p$ and $q$ be coprime integers 
(in particular, if one of them is $0$ then the other is $\pm1$).
Let $d$ be a positive integer.
We say that $L'$ is obtined from $L$ by a {\it $(dp,dq)$-cabling along $S$
with the core removed} (resp. {\it remained\,})
if $L' = (L\setminus S)\cup S'$ (resp. $L' = L\cup S'$) where 
$S'$ is a disjoint union of circles 
$S'=S_1\cup\dots\cup S_d\subset\partial T$ such 
that for each $j$ we have
$[S_j]=p[S]$ in $H_1(T)$ and $\lk(S_j,S')=q$.
A link is called an {\it iterated torus link} or a
{\it solvable link} if it can be successively
obtained from an unknot by these operations and maybe reversing the
orientations of some components.

To work with iterated toric links, we shall use the 
language of splice diagrams introduced by Eisenbud and Neumann in 
[\refEN]. A splice diagram is a tree $\Gamma$ which has no vertices
of valence $2$ and which is decorated as follows. Some of leaves 
(i.e. vertices of valence $1$) of $\Gamma$ are depicted as arrowheads,
a sign $\pm1$ is attrubuted to each arrowhead, 
and an integer $w(v,e)$ called the {\it edge weight} 
is attributed to each pair $(v,e)$ where $v$ is a {\it node} (i.e.
a vertex of valence $\ge3$) and $e$ is an edge incident to $v$.
Depicting splice diagrams,
we write $w(e,v)$ at the beginning of $e$ at $v$;
if the sign of an arrowhead is $+1$ then we do not write it.
The arrowhead vertices of the splice diagram of a link 
correspond to the link components.

We shall define the splice diagram of an iterated 
torus link inductively as follows.
The splice diagram of the unknot is 
$\,\bullet\!\!\!-\!\!\!-\!\!\!{\to}$\,.
Let $\boxed{\,\;\Gamma\;}\!{-}\!\!\!{-}\!\!{\overset{\eps}\to{\!\to}}$\,
with $\eps=\pm1$ be 
the splice diagram of $L$ and the depicted arrowhead corresponds to $S$.
If $L'$ is obtained from $L$ by a $(dp,dq)$-cabling along $S$ 
with the core removed (resp. remained) then the splice diagram of $L'$ is as
in Figure {\figCoreRemoved} (resp. in Figure \figCoreRemained).
Reversing of the orientation of a link component corresponds to
changing of the sign of the corresponding arrowhead.

\midinsert 
\epsfxsize 115mm
\centerline{\epsfbox{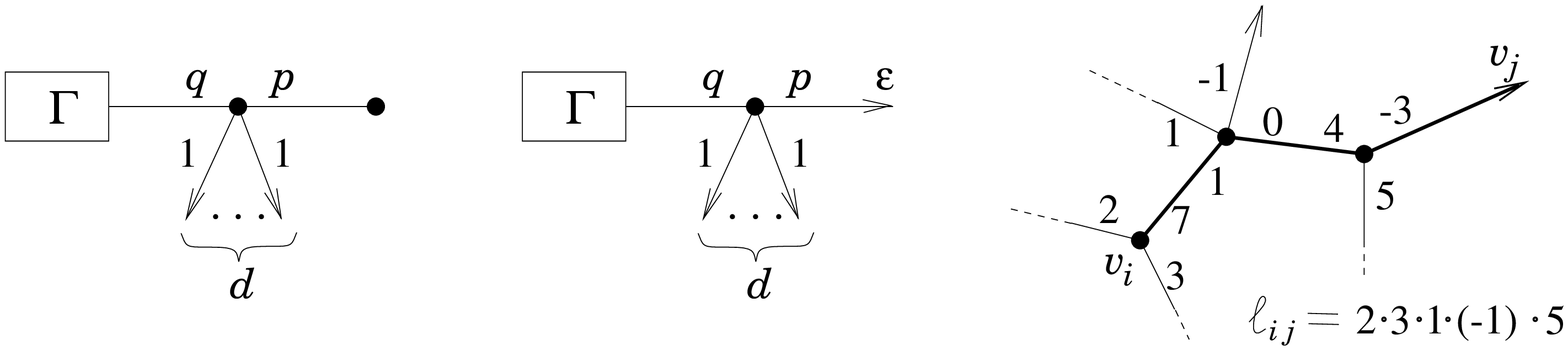}}
\botcaption{
		Fig.~\figCoreRemoved 
		\hbox to 25mm{}
		Fig.~\figCoreRemained
		\hbox to 27mm{}
		Fig.~\figlij\hbox to 5mm{}
          } 
\endcaption
\endinsert

A splice diagram is not determined by a link but Eisenbud and Neumann
[\refEN] described an equivalence relation between splice diagrams such that
the equivalence calss of a splice diagram is uniquely determined by a link.
Thay gave also a formula for the Alexander polynomial of an iterated torus link
% under the condition that it is fiberable 
in terms of its splice diagram
(in fact, this was done in [\refEN] for
graph links in homology spheres).
Neumann [\refNeumannJKTR] proved an analogous formula (only in the case of
fiberable links $S^3$) for the Conway polynomial.
%which is determined by the Alexander polynomial only up to a sign. 
To write this formula, we need the following notation.

Let $\Gamma$ be a splice diagram. Let us denote its vertices by
$v_1,\dots,v_n, v_{n+1},\dots,v_k$
with $v_1,\dots,v_n$ being arrowheads
(denote their signs by $\eps_1,\dots,\eps_n$)
and $v_{n+1},\dots,v_k$ being the remaining vertices
(denote their valences by $\delta_{n+1},\dots,\delta_k$).
For any two distinct vertices $v_i$ and $v_j$ of $\Gamma$ let us 
denote by $s_{ij}$ the simple path in $\Gamma$ joining $v_i$ to $v_j$,
including $v_i$ and $v_j$, and define (see Figure \figlij)
$$
	\ell_{ij}=\!\!\!\prod_{v\in s_{ij},e\not\subset s_{ij}}\!\!\!\!w(v,e),
	\quad i,j=1,\dots,k;
	\qquad
	m_i = \sum_{j=1}^n \ell_{ij}\eps_j, \quad i=n+1,\dots,k.
$$

% \example{ Example \exampleSD } If $\Gamma$ is the splice diagram
% depicted in Figure {\figSD} 

\proclaim{ Theorem \thEN } 
Let $L$ be a solvable link and
$\Gamma$ its splice diagram.
%$V$ the set of all non-arrowhead vertices of $\Gamma$.
\smallskip

a). {\rm[\refEN].} $L$ is fiberable if %and only if 
$m_i\ne 0$ for any $i=n+1,\dots,k$.
\smallskip

b). {\rm[\refNeumannJKTR].} If $L$ is fiberable then 
$$
	\Omega_L(t) = \eps_1\dots\eps_n \cdot (t-t^{-1}) \prod_{i=n+1}^k 
	\big(t^{m_i} - t^{-m_i}\big)^{\delta_i-2}
			%{\operatorname{valence}(v)-2}
							\eqno(\eqConway)
$$
\smallskip

\endproclaim

\proclaim{ Corollary \corEN }
Let $L$ be a solvable link and
$\Gamma$ its splice diagram.
Suppose that $m_i\ne0$ as soon as $\delta_i=1$
{\rm(}i.e. no denominator in the right hand side 
of {\rm(\eqConway)} vanishes{\rm)}. 
Then the equality {\rm(\eqConway)} holds.
\endproclaim

\demo{ Proof }
If $m_i\ne0$ for all $i=n+1,\dots,k$ then the 
result follows from Theorem \thEN.
If $m_i=0$ for some $i$ then $\Omega_L=0$ by Eisenbud-Neumann formula
for the multivariable Alexander polynomial [\refEN; Theorem 12.1].
\qed\enddemo

\remark{ Remark \remCimasoni }
Recently, Cimasoni [\refCimasoni] proved the analogue of Eisenbud-Neumann
formula for the multivariable Conway potential function for any 
graph link $L$ in a homology sphere:
$$
	\nabla_L(t_1,\dots,t_n) = \eps_1\dots\eps_n 
	\prod_{i=n+1}^k 
	\big(
	  t_1^{ \ell_{i1}\eps_1}\dots t_n^{ \ell_{in}\eps_n}
	- t_1^{-\ell_{i1}\eps_1}\dots t_n^{-\ell_{in}\eps_n}
	\big)^{\delta_i-2}
$$
(the terms $t_1^0\dots t_n^0 - t_1^0\dots t_n^0$ should be formally 
cancelled against each other before being set equal to zero).
Since $\Omega_L(t) = (t-t^{-1})\nabla_L(t,\dots,t)$, this formula 
allows one to compute $\Omega_L$ in the cases not covered
by Corollary {\corEN}, i.e. in the cases when (\eqConway) does not
make sense. 

For example (cp. [\refEN; p.~97]), if 
$\Gamma = $
$\smallmatrix -1\;\,1\;\;1\,\;\;\;\;1\\
\epsfxsize 13mm\epsfbox{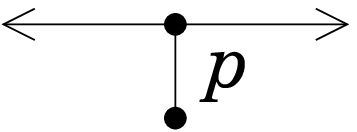}
\endsmallmatrix$
then $\nabla_L = -(u^p - u^{-p})/(u-u^{-1}) = -(u^{p-1}+\dots+u+1)$
where $u=t_1/t_2$. Hence, $\nabla_L(t,t)=-p$ and
$\Omega_L =(t^{-1}-t)p$ which well
agrees with the fact that a Seifert matrix of $L$ is the 
$1\times1$-matrix $(p)$.
%
%However, in this paper we use only Corollary \refEN.
\endremark

\midinsert 
\epsfxsize 125mm
\centerline{\epsfbox{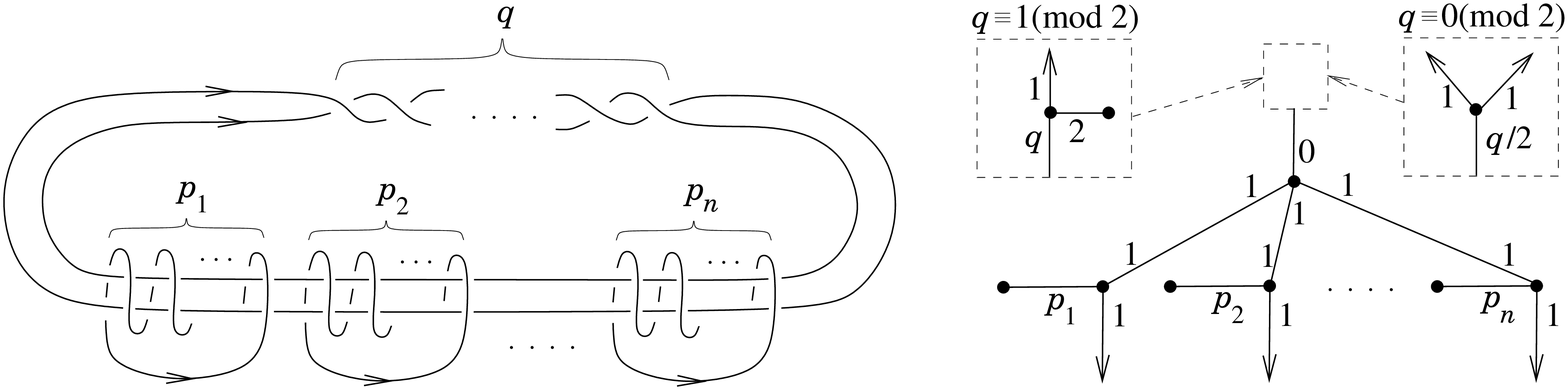}}
\botcaption{
              \hbox to 5mm{} Fig.~\figTI  
	      \hbox to 53mm{} Fig.~\figSD 
          } 
\endcaption
\endinsert

The following lemma provides an example of application of the above formulas.
It will be used in the next section.

\proclaim{ Lemma \lemTI} Let $L$ be the link in Fig.~\figTI\ where
$q,p_1,\dots,p_n\in\Z$
{\rm(}if some of the parameters $q,p_1,\dots,p_n$ are negative,
the corresponding positive crossings in Fig.~\figTI\ should
be replaced with negative ones;
if $p_j=0$ then the corresponding component is a simple closed curve
disjoint from the rest of the picture{\rm)}.
If $n>1$ then $\det_L=0$.
\endproclaim

\demo{ Proof } 
If $p_j=0$ for some $j$ then the link splits 
(has a component separated by a 2-sphere from the others) and
$\Omega_L=0$. 
Assume that $p_j\ne0$ for all $j$.
Then the splice diagram of $L$ is as in Fig.~\figSD.
%Set $m=2q+p_1+\dots+p_n$.
By Theorem \thEN, $L$ is fiberable unless 
$$
           q+p_1+\dots+p_n=0                   \eqno(\nonfib)
$$
and if it is fiberable then 
$\Omega_L(t) = (t-t^{-1})\prod_v\omega_v(t)$
where $v$ runs over the non-arrowhead vertices of valence $\ne2$.
Let $u$ be the vertex of valence $n+1$ (whose outcoming edges are weighted
by $0,1,\dots,1$). If $n>1$ then $\omega_u(t)=(t^2-t^{-2})^{n-1}$, hence,
$\omega_u(i)=0$. Thus, if $n>1$ and (\nonfib) does not hold then $\det_L=0$. 

In the case when (\nonfib) holds, set $L_0=L$ and 
let $L_\pm$ be the link obtained from
$L$ by replacing $q$ with $q\pm1$. Then (\nonfib) does not hold for 
$L_\pm$, and we have $\det L_\pm=0$. 
Hence, $\det L=0$ by (\skeindet).
\qed
\enddemo

% 
% \remark{ Remark } 
% 1. Since the sign of $\det L$ is not used, the
% proof of Lemma \lemTI\ does not really require the results of 
% [\refNeumannJKTR] but
% only the comutation of the Alexander polynomial from [\refEN].
% 
% 2. %It is not difficult to
% One can 
% prove Lemma \lemTI\ 
% directly by multiple induction using (\skeindet).
% \endremark
% 

%%%%%%%%%%%%%%%%%%%%%%%%%%%%%%%%%%%%%%%%%%%%%%%%%%%%%%%%%%%%%%%%%%%%%%%%%%%%
%%%%%%%%%%%%%%%%%%%%%%%%%%%%%%%%%%%%%%%%%%%%%%%%%%%%%%%%%%%%%%%%%%%%%%%%%%%%
%%%%%%%%%%%%%%%%%%%%%%%%%%%%%%%%%%%%%%%%%%%%%%%%%%%%%%%%%%%%%%%%%%%%%%%%%%%%
%%%%%%%%%%%%%%%%%%%%%%%%%%%%%%%%%%%%%%%%%%%%%%%%%%%%%%%%%%%%%%%%%%%%%%%%%%%%
%%%%%%%%%%%%%%%%%%%%%%%%%%%%%%%%%%%%%%%%%%%%%%%%%%%%%%%%%%%%%%%%%%%%%%%%%%%%
%%%%%%%%%%%%%%%%%%%%%%%%%%%%%%%%%%%%%%%%%%%%%%%%%%%%%%%%%%%%%%%%%%%%%%%%%%%%

\head \sectMainLemma. Link theoretical lemma
\endhead

Let $S$ be an oriented unknotted circle in $S^3$ and
$L$ an oriented link which contains $S$ as a component.
Suppose that $L'$ is a link whose diagram is obtained from 
a diagram of $L$ by simultaneous replacing of a negative
crossing with a positive one and a positive crossing 
with a negative one where each of the two crossings
involves a segment of $S$ and a segment of $L\setminus S$
(see Fig.~\figLemmaA\ where $S$ is thicker than $L\setminus S$).

\midinsert 
\epsfxsize 120mm
\centerline{\epsfbox{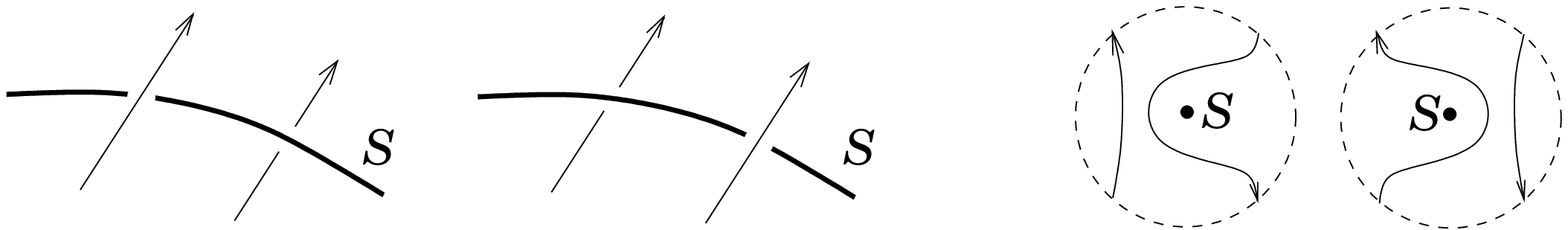}}
\centerline{\hbox to 5mm{} $L$ \hbox to 30mm{} $L'$ 
            \hbox to 36mm{} $L$ \hbox to 15mm{} $L'$
	   }
\botcaption{
              \hbox to 10mm{} Fig.~\figLemmaA  
	      \hbox to 53mm{} Fig.~\figLemmaB 
          } 
\endcaption
\endinsert

In other words, 
$L'$ is obtained from $L$ as follows.
Choose disjont embedded disks $D_1$ and $D_2$
such that
\roster
\item $a_i:=D_i\cap(L\setminus S)$ is a segment of $\partial D_i$, $i=1,2$;

\item each $D_i$ meets $S$ transversally at a single point;

\item the signs of the intersections are opposite 
(i.e. $D_1\cdot S=-D_2\cdot S$)
where the orientation of $D_i$ is chosen so that the orientation of 
$a_i$ induced from $L$ coincides with that induced from $\partial D_i$.
\endroster
Then $L'=(L\setminus(a_1\cup a_2))\cap(a'_1\cup a'_2)$
where $a'_i=\overline{\partial D_i\setminus a_i}$.

\proclaim{ Lemma \MainLemma } 
Let $q$ be any integer and $L_{2,q}$ (resp. $L'_{2,q}$)
be the $(2,q)$-cabling of $L$ (resp. of $L'$) along $S$ 
with the core removed. Then:

a). $\det L'_{2,q} = \det L_{2,q}$;

b). If $\det L_{2,q}\ne 0$ then $\Sign L'_{2,q} = \Sign L_{2,q}$.
\endproclaim

\demo{ Proof } a).
Let us identify $S^3$ with $\R^3\cup\{\infty\}$ so that $S$ is a line
and the projection of $L\setminus S$ onto a plane along the line $S$ is
non-degenerate.
One can assume that there exists a disk $D\subset\R^2$ 
centered at the projection of $S$ such that
the diagrams of $L\setminus S$ and $L'\setminus S$ 
(with respect to the chosen projection) coincide outside $D$ 
and have the form shown in Fig.~\figLemmaB\ inside $D$.

Let us use the induction by the number of crossings of the diagram of
$L\setminus S$.
Suppose that there are no crossings. If $L\setminus S$ is connected
then $L'=L$. Otherwise, each of $L_{2,q}$ and $L'_{2,q}$ has the form 
shown in Fig.~\figTI\ with $p_j=0,\pm1$ and $n>1$, 
hence, $\det{L_{2,q}}=\det{L'_{2,q}}=0$ by Lemma \lemTI.

Now suppose that the number of crossings of the diagram of
$L\setminus S$ is $N$ and suppose that the statement of the lemma is proved 
for all diagrams with less than $N$ crossings. 
Denote by $\pi$  the projection $\R^3\to \R^2$ (recall that $\pi(S)$ is
a single point).
Let $L_1,\dots,L_n$ be the copmonents of $L\setminus S$. 
On each $L_j$, let us choose
a point $x_j$ such that $\pi(x_j)$ lies on the boundary of the unbounded
component of $\R^2\setminus\pi(L_j)$ and $\pi(x_j)$ 
is not a double point of $\pi(L\setminus S)$.
Let $f:(L\setminus S)\setminus\{x_1,\dots,x_n\}\to\R$ maps 
homeomorphically each $L_j\setminus x_j$ onto the interval $(j-1,j)$.
Let $x,y\in L$ be such that $\pi(x)=\pi(y)$ and $f(x)<f(y)$. 
Say that the corresponding crossing of the diagram of $L\setminus S$ 
is {\it monotone} if $y$ is higher than $x$.

Consider the family 
of all diagrams which differ from the diagram of
$L$ by changing the signs of crossings 
(the number of such diagrams is $2^N$).
Let us prove the statement of the lemma for the diagrams from this family
using the induction by the number of non-monotone crossings
(with respect to the same choice of $x_1,\dots,x_n$).

{\it Base case of the induction.} 
Suppose that all the crossings of $L$ are monotone.
Then each of $L_{2,q}$ and $L'_{2,q}$ has the form shown in Fig.~\figTI\
where $p_j$ is the linking number $\lk(L_j,S)$
(if $p_j=0$ then $L_j$ should be shown in Fig.~\figTI\ by a simple curve
disjoint from the rest of the picture).
Hence, $L=L'$ for $n=1$ and, by Lemma \lemTI, we have 
$\det L_{2,q}=\det L'_{2,q}=0$ for $n>1$.

{\it Step of the induction.}
Consider a non-monotone crossing of the diagram of $L\setminus S$.
Let it be positive (the case of a negative non-monotone crossing
can be treated similarly).
Set $L_+=L$, $L'_+=L'$ and let $L_0$ and $L_-$ (resp. $L'_0$ and $L'_-$)
be obtained from $L_+$ (resp. $L'_+$) by 
changing the crossing according to Fig.~\figSkein.
By the induction hypothesis, we have 
$\det (L_0)_{2,q}=\det (L'_0)_{2,q}$
(because the number of crossings of $L_0$ is less than that of $L_+$) and
$\det (L_-)_{2,q}=\det (L'_-)_{2,q}$
(because the number of non-monotone crossings of $L_-$ 
is less than that of $L_+$). Hence, 
$\det (L_+)_{2,q}=\det (L'_+)_{2,q}$
by (\skeindet).

b). Follows from (a) and Corollary \corSkeinSigma.
\qed
\enddemo

%%%%%%%%%%%%%%%%%%%%%%%%%%%%%%%%%%%%%%%%%%%%%%%%%%%%%%%%%%%%%%%
%%%%%%%%%%%%%%%%%%%%%%%%%%%%%%%%%%%%%%%%%%%%%%%%%%%%%%%%%%%%%%%
%%%%%%%%%%%%%%%%%%%%%%%%%%%%%%%%%%%%%%%%%%%%%%%%%%%%%%%%%%%%%%%
%%%%%%%%%%%%%%%%%%%%%%%%%%%%%%%%%%%%%%%%%%%%%%%%%%%%%%%%%%%%%%%
%%%%%%%%%%%%%%%%%%%%%%%%%%%%%%%%%%%%%%%%%%%%%%%%%%%%%%%%%%%%%%%

\head \sectDet. Computation of some braid determinants
\endhead

\definition{ Notation \notBraid }
Let us denote the standard generators of the braid group $B_m$
by $\sigma_1,\sigma_2,\dots,\sigma_{m-1}$
and let us set:
$$
	\pi_{k,l} = \cases \sigma_k\sigma_{k+1}\dots\sigma_l, &k<l,\\
			\sigma_k\sigma_{k-1}\dots\sigma_l, &k>l,\\
			\sigma_k, &k=l,
		\endcases
\qquad\qquad
     \tau_{k,l}=\cases \pi_{l,k+1}^{-1}\pi_{k,l-1},  &k<l, \\ 
                       \pi_{l,k-1}^{-1}\pi_{k,l+1},       &k>l, \\ 
                       1,                                  &k=l, 
                \endcases 
$$
$\Delta_k = \pi_{1,k-1}\pi_{1,k-2}\dots\pi_{1,2}\sigma_1$.
\enddefinition

For positive interegs $n$, $J$, $k$ such that $n\equiv J \mod 2$
and non-negative integers $\alpha_1,\dots,\alpha_J$
let us set $m=2k+1$ and define the braid 
$b=b_{n,k}^J(\alpha_1,\dots,\alpha_J)\in B_m$
as
$$
    b=b_1\dots b_J\,\Delta^n, \qquad
    b_j = \cases
    		\sigma_{k-1}^{-\alpha_j}\tau_{k-1,k},
			&\text{ if $j$ is odd,}\\
    		\sigma_{k  }^{-\alpha_j}\tau_{k,k-1},
			&\text{ if $j$ is even.}
	  \endcases
$$
If, moreover, $k\ge2$, we define the braid 
$c=c_{n,k}^J(\alpha_1,\dots,\alpha_J)\in B_{2k+1}$
as
$$
    c=c_1\dots c_J\,\Delta^n, \qquad
    c_j = \cases
    		\sigma_{k-2}^{-\alpha_j}\tau_{k-2,k+1},
			&\text{ if $j$ is odd,}\\
    		\sigma_{k+1}^{-\alpha_j}\tau_{k+1,k-2},
			&\text{ if $j$ is even.}
	  \endcases
$$

For $\vec\alpha=(\alpha_1,\dots,\alpha_J)$, $\alpha=\alpha_1+\dots+\alpha_J$,
let us set
$$
	\tilde b_{n,k}^J(\vec\alpha)=i^\alpha\det b_{n,k}^J(\vec\alpha),
	\qquad
	\tilde c_{n,k}^J(\vec\alpha)=i^\alpha\det c_{n,k}^J(\vec\alpha),
                                                      \eqno(\defdnkJ)
$$

\proclaim{ Lemma \lemDet }
a). Let $n,k,J\ge1$ and $n\equiv J \mod 2$. Then
$$
   \det b_{n,k}^J(1,\dots,1) = 
	\cases
		4	&\text{if $n\equiv J+2\equiv0\mod4$ and $k=1$,}\\
		4^k	&\text{if $n+2\equiv J\equiv0\mod4$,}\\
		%\Big((-1)^{(n+1)/2}\cdot 2i\Big)^{k+1}
		(-2i^n)^{k+1}
			&\text{if $n$ is odd and $J\equiv n+2k\mod4$,}\\
		0	&\text{otherwise.}
	\endcases
$$

b). Let $n,J\ge1$, $k\ge2$, and $n\equiv J \mod 2$. Then
$$
   \det c_{n,k}^J(1,\dots,1) = 
	\cases
		4^k&\text{if $n+2\equiv J\equiv0\mod4$}\\
		%\Big((-1)^{(n+1)/2}\cdot 2i\Big)^{k+1}
		-(-2i^n)^{k+1}
			&\text{if $n$ is odd and $J\equiv n+2k+2\mod4$,}\\
		0	&\text{otherwise.}
	\endcases
$$
\endproclaim

% \midinsert 
% \epsfxsize 125mm
% \centerline{\epsfbox{dnfb111.eps}}
% \botcaption{
%               \hbox to  5mm{} Fig.~\figBraid  
% 	      \hbox to 53mm{} Fig.~\figBraidSD 
%           } 
% \endcaption
% \endinsert

\demo{ Proof } a).
Let us show that the closure $L$ of $b_{n,k}^J$ is an iterated torus link
and its splice diagram is shown in Fig.~\figBraidSD.
Indeed, let us apply the procedure form Section \sectSplice.
We start with the $k$-th string and we add the strings
number $1,\dots,k-1$, $k+3,\dots,m$ as the $(2,n)$-  
or $(1,n/2)$- (according to the parity of $n$) -cables
along it (see Fig.~\figBraid).
Then we add the $(k+1)$-th string as the 
$(1,(n-J)/2)$-cable along the $k$-th string and finally, we add the 
$(k+2)$-th string as the $(1,(n+J)/2)$-cable along the $(k+1)$-th string.

\midinsert 
\epsfxsize 125mm
\centerline{\epsfbox{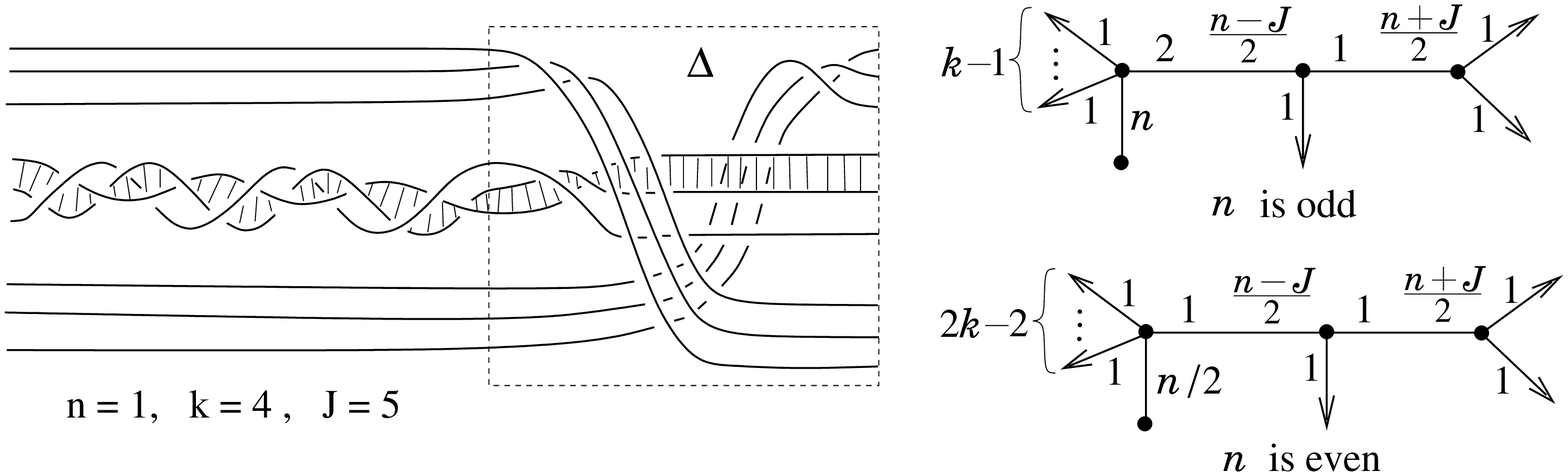}}
\botcaption{
              \hbox to  5mm{} Fig.~\figBraid. 
	      The braid $b_{n,k}^J$ %for $n=1$, $k=4$, $J=5$
	      \hbox to 22mm{} Fig.~\figBraidSD. The s.-d. of 
	      $b_{n,k}^J\!\!\!\!\!\!\!\!\!$
          } 
\endcaption
\endinsert

By Corollary \corEN, we have 
$$
	\Omega_L(t) = {
	  (t-t^{-1})\cdot\omega_1^{k-1}
	  \cdot(t^{m_2}-t^{-m_2})(t^{m_3}-t^{-m_3})
	  \over
	  t^{m}-t^{-m}
	}
	 					\eqno(\eqOmegaL)
$$
where
$$
	\omega_1={(t^{nm}-t^{-nm}) }
	\;\text{ if $n$ is odd,}
	\quad
	\omega_1={(t^{nm/2}-t^{-nm/2})^2 }
	\;\text{ if $n$ is even}.
$$
and 

%\smallskip
%$\qquad\qquad m=2k+1$, 

\smallskip
$\qquad\qquad m_2=n(k-1)+{3\over2}(n-J)={mn\over2}-{3\over2}J$,

\smallskip
$\qquad\qquad m_3=n(k-1)+{1\over2}(n-J)+(n+J)={mn\over2}+{1\over2}J$.

%
%\smallskip
%{Case 1.} $n\equiv0\mod4$. If $k>1$ then $\omega_1(i)=0$.
%If $k=1$ then 
%
%\smallskip
%{Case 2.} $n$ is odd. Since $m_3 - m_2 = 2J$
%
\medskip

b). The proof is similar to that of the part (a).
The closure $L$ of $c_{n,k}^J$ is an iterated torus link
and its splice diagram is shown in Fig.~\figBraidSDc. 
Indeed, if we remove the $(k+2)$-th and $(k+3)$-th strings,
we obtain the braid $b_{n,k-2}^2$ whose splice diagram is
computed above and these two strings can be considered as the
$(2,n)$- or $(1,n/2)$-cables along the $(k+1)$-th string.

\midinsert 
\epsfxsize 125mm
\centerline{\epsfbox{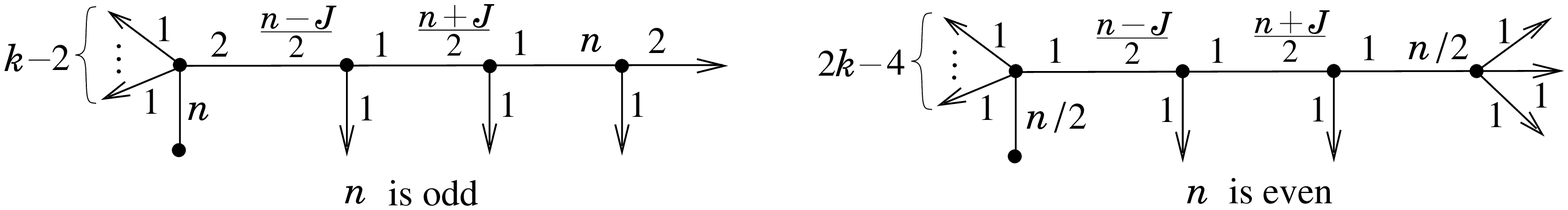}}
\botcaption{
	      Fig.~\figBraidSDc. The splice diagram of $c_{n,k}^J$
          } 
\endcaption
\endinsert

Thus, by Corollary \corEN, we obtain
$$
	\Omega_L(t) = {
	  (t-t^{-1})\cdot\omega_1^{k-2}\cdot
          (t^{m_2}-t^{-m_2})(t^{m_3}-t^{-m_3})\cdot\omega_4
	  \over
	  t^{m}-t^{-m}
	}
$$
where 
$$
\xalignat3
	&\omega_1={(t^{nm}-t^{-nm})},&
	&\omega_4=t^{2m_4}-t^{-2m_4} &
	&\text{ if $n$ is odd, }
	\\
	&\omega_1={(t^{nm/2}-t^{-nm/2})^2},&
	&\omega_4=(t^{m_4}-t^{-m_4})^2 &
	&\text{ if $n$ is even. }
\endxalignat
$$
and 

%\smallskip
%$\qquad \qquad m=2k+1$, 

\smallskip
$\qquad \qquad m_2=n(k-2)+{5\over2}(n-J)={mn\over2}-{5\over2}J$, 

\smallskip
$\qquad \qquad m_3=n(k-2)+{1\over2}(n-J)+2(n+J)={mn\over2}+{3\over2}J$, 

\smallskip
$\qquad \qquad m_4=n(k-2)+{1\over2}(n-J)+{1\over2}(n+J)+{3\over2}n={mn\over2}$

\medskip\noindent
(we see a posteriori that $\omega_4=\omega_1$).
\qed
\enddemo

\proclaim{ Lemma \lemDetBase }
a). If $n$ is odd and $k\ge 1$ then
$
	\det b_{n,k}^1(0) = (2i^n)^k. % = \Big((-1)^{(n-1)/2}\cdot 2i\Big)^k.
$

b). If $n$ is even and $k>0$ then 
$$
	\det b_{n,k}^2(0,0) = %%\det c_{n,k}^2(0,0) = 
	i\,\det b_{n,k}^2(1,0) = (2-2i^n)^k =
	\cases 
		4^k, 	&n\equiv2\mod4,\\ 
		0, 	&n\equiv0\mod4,
	\endcases
$$

c). If $k\ge2$ then the braids 
$c_{n,k}^1(0)$, $c_{n,k}^2(0,0)$, and $c_{n,k}^2(1,0)$ 
are conjugate to 
$b_{n,k}^1(0)$, $b_{n,k}^2(0,0)$, and $b_{n,k}^2(1,0)$ 
respectively. 
\endproclaim

\demo{ Proof }
a). Since $b_{n,k}^1(0) = \sigma_k^{-1}\sigma_{k-1}\Delta^n = 
\sigma_k^{-1} \Delta^n \sigma_k$ is conjugate to $\Delta^n$,
it represents the toric link whose diagram is
depicted in Figure \figBaseSD. Hence,
$$
	\Omega_{b_{n,k}^1(0)}(t) = (t-t^{-1})\cdot
	{(t^{m n}-t^{-m n})^k\over t^{m}-t^{-m}},
\qquad\text{where $m=2k+1$}.
$$

\midinsert 
\epsfxsize 120mm
\centerline{\epsfbox{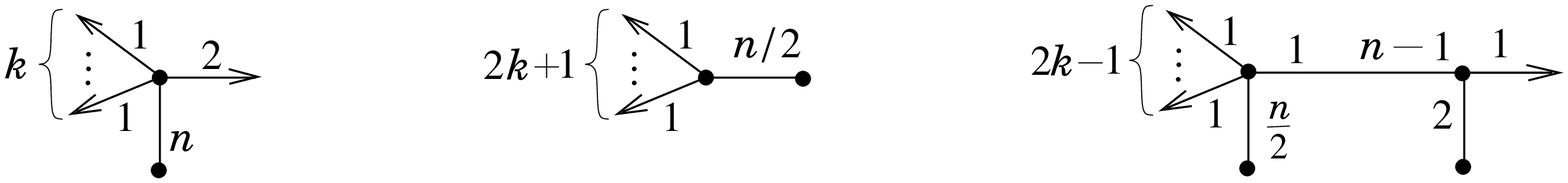}}
\noindent
	$b_{n,k}^1(0)=\sigma_k^{-1}\Delta^n\sigma_k$
\hskip 10mm 
	$b_{n,k}^2(0,0)=\Delta^n$
\hskip 30mm 
	$b_{n,k}^2(1,0)$
\botcaption{ Fig.~\figBaseSD  } 
\endcaption
\endinsert

b). The links represented by the braids 
$b_{n,k}^2(0,0)%=c_{n,k}^2(0,0)
=\Delta^n$ and 
$b_{n,k}^2(1,0)$ are solvable. Their splice diagrams are as in
Figure \figBaseSD, so we have
$$
\split
	&\Omega_{b_{n,k}^2(0,0)}(t) = (t-t^{-1})\cdot
	{(t^{m n/2}-t^{-m n/2})^{2k}\over t^{m}-t^{-m}},
\\
	&\Omega_{b_{n,k}^2(1,0)}(t) = (t-t^{-1})\cdot
	{(t^{m n/2}-t^{-m n/2})^{2k-1}\over t^{m}-t^{-m}}\cdot
	{t^{2m_2}-t^{-2m_2}\over t^{m_2}-t^{-m_2}}
\endsplit
$$
where
%$m=2k+1$, 
$m_2=(2k-1){n\over2}+(n-1) = {mn\over2}-1>0$.

c). Easy to check.
\qed
\enddemo

%%%%%%%%%%%%%%%%%%%%%%%%%%%%%%%%%%%%%%%%%%%%%%%%%%%%%%%%%%%%%%%%%%%%%
%%%%%%%%%%%%%%%%%%%%%%%%%%%%%%%%%%%%%%%%%%%%%%%%%%%%%%%%%%%%%%%%%%%%%
%%%%%%%%%%%%%%%%%%%%%%%%%%%%%%%%%%%%%%%%%%%%%%%%%%%%%%%%%%%%%%%%%%%%%
%%%%%%%%%%%%%%%%%%%%%%%%%%%%%%%%%%%%%%%%%%%%%%%%%%%%%%%%%%%%%%%%%%%%%
%%%%%%%%%%%%%%%%%%%%%%%%%%%%%%%%%%%%%%%%%%%%%%%%%%%%%%%%%%%%%%%%%%%%%
%%%%%%%%%%%%%%%%%%%%%%%%%%%%%%%%%%%%%%%%%%%%%%%%%%%%%%%%%%%%%%%%%%%%%
%%%%%%%%%%%%%%%%%%%%%%%%%%%%%%%%%%%%%%%%%%%%%%%%%%%%%%%%%%%%%%%%%%%%%

\head \sectCycPol. Skein systems of cyclically symmetric polynomials
\endhead

\definition{ Definition \defCycPol }
Let $\nu=1$ or $2$.
A sequence of polynomials $\{f_J(x_1,\dots,x_J)\}$, 
${J=\nu,\nu+2,\nu+4,\dots}$,
is called a {\it skein system of cyclically symmetric polynomials 
of the parity}
$\nu$ (or, just a {\it skein system of parity $\nu$}) 
if the following conditions hold
\roster
\item "({\it i}\,)"
	$\deg_{x_1}f_J \le 1$;
\item "({\it ii}\,)"
	$f_J(x_1,\dots,x_J) = f_J(x_J,x_1,\dots,x_{J-1})$;
\item "({\it iii}\,)"
	$f_J(x_1,0,x_3,\dots,x_J) = f_{J-2}(x_1+x_3,x_4,\dots,x_J)$
       for $J\ge\nu+2$.
\endroster
% A skein system $\{f_J(x_1,\dots,x_J)\}$ is called {\it $q$-periodic}
% if $f_{J+q}(1,\dots,1) = f_J(1,\dots,1)$.
\enddefinition

\proclaim{ Lemma \lemDetCP } Let $f_J=\tilde b_{n,k}^J$ or $\tilde c_{n,k}^J$
{\rm(}see Section \sectDet{\rm)}.

a). If we fix positive integers $n,J,k$ such that $n\equiv J \mod 2$  
then $f_J(\alpha_1,\dots,\alpha_J)$ are values at integral points
of a polynomial
{\rm(}we denote it also $f_J${\rm)}. 

b). For fixed $n,k>0$, the sequences $\{f_J\}_{J>0,\,J\equiv n(2)}$ 
is a skein system of the parity $n$.
\endproclaim

\demo{ Proof }
When all $x_j$ are integer, ({\it ii}\,) and ({\it iii}\,) are
evident and it follows from (\skeindet) that 
$f_J$ is linear with respect to each variable.
\qed
\enddemo

\proclaim{ Lemma \lemCycPolA }
Let $\calF=\{f_J\}$ be a skein system of parity $\nu$.
Let $c_J = f_J(1,\dots,1)$.

If $\nu=1$ then $\calF$ is uniquely determined by the sequence
$\{c_1,c_3,c_5,\dots\}$ and the number $c_0=f_1(0)$.

If $\nu=2$ then $\calF$ is uniquely determined by the sequence
$\{c_2,c_4,c_6\dots\}$ and the numbers $c_0=f_2(0,0)$, $c_1=f_2(1,0)$.
\endproclaim

\demo{ Proof } 
Suppose we proved the uniqueness of $f_{J-2}$.
By ({\it i\,}) -- ({\it ii\,}), $f_J$ is determined by its values at
the vertices of the unit cube. The value at $(1,\dots,1)$ is $c_J$. 
The values at other vertices can be expressed by 
({\it ii\,}) and ({\it iii\,}) in terms of $f_{J-2}$.
\qed
\enddemo

Let us define symmetric $(J\times J)$-matrices 
$A_J^{\pm}(x_1,\dots,x_J)$ as follows.
Let $E_{i,j}$ be the ($J\times J$)-matrix whose $(k,l)$-th entry is 
$\delta_{ki}\delta_{jl}$.
Set
$$
	A_J^{\pm}(x_1,\dots,x_J) 
		= -2\,\left(\sum_{i=1}^J x_i E_{i,i}\right) +
			\left(\sum_{i=1}^{J-1} E_{i,i+1}+E_{i+1,i}\right)
			\pm (E_{1,J} + E_{J,1}).
$$
Thus,
$$
\split
	&A_1^\pm = (\pm2-2x_1),\; 
	A_2^+ = \left(\matrix -2x_1&2\\2&-2x_2 \endmatrix\right),\;
	A_2^- = \left(\matrix -2x_1&0\\0&-2x_2 \endmatrix\right),
\\
	&A_3^\pm = \left(\matrix -2x_1&1&\pm1\\
				1&-2x_2&1\\
				\pm1&1&-2x_3 \endmatrix\right),\;
	A_4^\pm = \left(\matrix -2x_1&1&0&\pm1\\
	1&-2x_2&1&0\\0&1&-2x_3&1\\
	\pm1&0&1&-2x_4 \endmatrix\right),\;\dots
\endsplit
$$

\proclaim{ Lemma \lemDetB} 
$\det A_J^-(x_1,\dots,x_J)=\det A_J^+(x_1,\dots,x_J)+(-1)^J\cdot 4$.
\endproclaim

\demo{ Proof } 
The determinant of the matrix obtained from $A_J^+$
by replacing the $(1,J)$-th and $(J,1)$-th entry with $u$, 
is a quadratic function of $u$ whose linear term is equal to
$(-1)^{J+1}\cdot2u$.
\qed
\enddemo

\proclaim{ Lemma \lemDetA} 
	$\det A_J^\pm(x_1,0,x_3,\dots,x_J) 
	= -\det A_{J-2}^\mp(x_1+x_3,x_4,\dots,x_J)$.
\endproclaim

\demo{ Proof } %Follows from Lemma \lemDetB\ and the fact that
Indeed, we have $\det A_J^\pm(x_1,0,x_3,\dots,x_J)\,=$
\smallskip
$$
\split 
&
 =\left|\matrix     -2x_1 &   1  &  0   & \dots & \pm1\\
                       1  &   0  &  1   & \dots &   0 \\
		       0  &   1  &-2x_3 & \dots &   0 \\
		    \vdots&\vdots&\vdots&       &     \\
		     \pm1 &   0  &  0   &       &   
		                                    \endmatrix\right|
  =\left|\matrix       0  &   1  &  x_1 & \dots & \pm1\\
                       1  &   0  &  1   & \dots &   0 \\
		      x_1 &   1  &-2x_3 & \dots &   0 \\
		    \vdots&\vdots&\vdots&       &     \\
		     \pm1 &   0  &  0   &       &
		                                    \endmatrix\right| =
\\
& =\left|\matrix       0  &   1  &  x_1     & \dots & \pm1\\
                       1  &   0  &  0       & \dots &   0 \\
		      x_1 &   0 &-2(x_3+x_1)& \dots & \mp1\\
		    \vdots&\vdots&\vdots    &       &     \\
		     \pm1 &   0  & \mp1     &       &
		                                    \endmatrix\right|
  =\left|\matrix0&1\\1&0\endmatrix\right|\cdot
   \det A_{J-2}^\mp(x_1+x_3,x_4,\dots,x_J).\qed
\endsplit
$$
\enddemo

Let us denote 
$
	a_J^\pm=\cases\pm\det A_J^\pm &\text{ if $J\equiv 0$ or $1\mod 4$},\\
		      \mp\det A_J^\mp &\text{ if $J\equiv 2$ or $3\mod 4$}.
		\endcases
$
%
%  a_1^+(0)     =  det A_1^+(0) = 2
%  a_1^+(1)     =  det A_1^+(1) = 0 
%  a_3^+(1,1,1) = -det A_3^-(1,1,1) = -4 (-1)^3 = 4
%             . . . 
%
%  a_1^-(0)     = -det A_1^-(0) = 2
%  a_1^-(1)     = -det A_1^-(1) = -4 (-1)^1 = 4
%  a_3^-(1,1,1) =  det A_3^+ = 0
%             . . . 
%
%  a_2^+(x,y)     = -det A_2^-(x,y) = -(4xy)
%  a_2^+(0,0)     = 0
%  a_2^+(1,0)     = 0
%  a_2^+(1,1)     = -det A_2^-(1,1) = -4 (-1)^2 = -4
%  a_4^+(1,1,1,1) =  det A_4^+(1,1,1,1) = 0
%
%  a_2^-(x,y)     =  det A_2^+(x,y) = (4xy-4)
%  a_2^-(0,0)     = -4
%  a_2^-(1,0)     = -4
%  a_2^-(1,1)     =  det A_2^+(1,1) = 0
%  a_4^-(1,1,1,1) = -det A_2^-(1,1,1,1) = -4 (-1)^4 = -4
%

\proclaim{ Corollary \corCycPol }

\noindent
$\{a_1^+,a_3^+,a_5^+\dots\}$ and $\{a_1^-,a_3^-,a_5^-\dots\}$ are odd skein systems;
%$\calA_0^+$ and $\calA_0^-$ are even skein systems;

\noindent
$\{a_2^+,a_4^+,a_6^+\dots\}$ and $\{a_2^-,a_4^-,a_6^-\dots\}$ are even skein systems;
%$\calA_1^+$ and $\calA_1^-$ are odd skein systems; 
\qed\endproclaim

\proclaim{ Lemma \lemSignDetA }
	a). $\det A_J^+(1,\dots,1)=0$;
\smallskip

	b). $\det A_J^-(1,\dots,1)=(-1)^J\cdot4$.
\smallskip

	c). If $x_j\ge 1,\dots,x_J\ge 1$ and 
	$(x_1,\dots,x_J)\ne(1,\dots,1)$ then 
$$
	\sign\det A_J^+(x_1,\dots,x_J)=(-1)^J.
$$
	d). If $x_j\ge 1,\dots,x_J\ge 1$ then 
		$\;\sign\det A_J^-(x_1,\dots,x_J)=(-1)^J$.
\endproclaim

\demo{ Proof } a). $\det A_J^+=0$ because the sum of the rows is zero;

\smallskip
b). Follows from (a) and Lemma \lemDetB.

\smallskip
c). Set $A_J^+(x_1,\dots,x_J)=A+D$ where $A=A_J^+(1,\dots,1)$ and
$D=\diag(2-2x_1,\dots,2-2x_J)$.
Let $e_1,\dots,e_J$ be the standard base of $\R^J$. Then
$A$ defines a quadratic form on $\R^J$ whose restriction onto
$\langle e_1,\dots,e_{J-1}\rangle$ is negative definite 
(this is the Cartan matrix of the type $A$)
and the kernel of the form $A$ 
is generated by the vector $v=e_1+\dots+e_J$. The diagonal form
$D$ is non-positive and it is negative on $v$.
Thus, the form $A+D$ is negative definite.

\smallskip
d). Follows from (b) and the fact that the principal
$(J-1)\times(J-1)$-minors are negative definite (see the proof of (c)). 
\qed\enddemo

%
%  a_1^+(0)     =  det A_1^+(0) = 2
%  a_1^+(1)     =  det A_1^+(1) = 0 
%  a_3^+(1,1,1) = -det A_3^-(1,1,1) = -4 (-1)^3 = 4
%  a_5^+(1...1) =  det A_5^+(1...1) = 0
%  a_7^+(1...1) = -det A_7^-(1...1) = -4 (-1)^7 = 4
%  a_9^+(1...1) =  det A_9^+(1...1) = 0
%             . . . 
%
%  a_1^-(0)     = -det A_1^-(0) = 2
%  a_1^-(1)     = -det A_1^-(1) = -4 (-1)^1 = 4
%  a_3^-(1,1,1) =  det A_3^+(1,1,1) = 0
%  a_5^-(1...1) = -det A_5^-(1...1) = -4 (-1)^5 = 4
%  a_7^-(1...1) =  det A_3^+(1...1) = 0
%  a_9^-(1...1) = -det A_9^-(1...1) = -4 (-1)^9 = 4
%             . . . 
%
%  a_2^+(x,y)     = -det A_2^-(x,y) = -(4xy)
%  a_2^+(0,0)     = 0
%  a_2^+(1,0)     = 0
%  a_2^+(1,1)     = -det A_2^-(1,1) = -4 (-1)^2 = -4
%  a_4^+(1,1,1,1) =  det A_4^+(1,1,1,1) = 0
%  a_6^+(1,...,1) = -det A_6^-(1,...,1) = -4 (-1)^6 = -4
%  a_8^+(1,...,1) =  det A_8^+(1,...,1) = 0
%             . . . 
%
%  a_2^-(x,y)     =  det A_2^+(x,y) = (4xy-4)
%  a_2^-(0,0)     = -4
%  a_2^-(1,0)     = -4
%  a_2^-(1,1)     =  det A_2^+(1,1) = 0
%  a_4^-(1,1,1,1) = -det A_4^-(1,1,1,1) = -4 (-1)^4 = -4
%  a_6^-(1,...,1) =  det A_6^+(1,...,1) = 0
%  a_8^-(1,...,1) = -det A_4^-(1,...,1) = -4 (-1)^8 = -4
%             . . . 
%

\proclaim{ Corollary \corBaseCycPol } The skein systems from 
Corollary {\corCycPol} satisfy the following initial conditions.
$$
\xalignat3
	&a_1^+(0)=2,\qquad
	&&a_J^+(1,\dots,1)=2(i^{J+1}+1),&& J=1,3,5,\dots
	\\
	&a_1^-(0)=2,\qquad
	&&a_J^-(1,\dots,1)=2(i^{J-1}+1),&& J=1,3,5,\dots
	\\
	&a_2^+(0,0)=a_2^+(1,0)=0,\qquad
	&&a_J^+(1,\dots,1)=2(i^J-1),&& J=2,4,6,\dots
	\\
	&a_2^-(0,0)=a_2^-(1,0)=-4,\qquad
	&&a_J^-(1,\dots,1)=-2(i^J+1),&& J=2,4,6,\dots\qed
\endxalignat
$$

\endproclaim

\proclaim{ Corollary \corDnkJ } Let $n,k,J$ be positive integers and let
$\vec\alpha=(\alpha_1,\dots,\alpha_J)$, $\alpha_j\ge 0$.

\smallskip
\noindent
If $n\equiv 0\mod4$ and $k>1$ then 
$$
\tilde b_{n,k}^J(\vec\alpha)=
\tilde c_{n,k}^J(\vec\alpha)=0.
$$
If $n\equiv 0\mod4$ and $k=1$ then
$$
\tilde b_{n,k}^J(\vec\alpha) = 
a_J^+(\vec\alpha).
%\cases
%	+\det A_J^+ &\text{if $J\equiv 0\mod4$,}\\
%	-\det A_J^- &\text{if $J\equiv 2\mod4$.}
%\endcases
$$
If $n\equiv 2\mod4$ then 
$$
\tilde b_{n,k}^J(\vec\alpha) =
\tilde c_{n,k}^J(\vec\alpha) =
 -4^{k-1}a_J^-(\vec\alpha).
%\cases
%	+4^{k-1}\det A_J^- &\text{if $J\equiv 0\mod4$,}\\
%	-4^{k-1}\det A_J^+ &\text{if $J\equiv 2\mod4$.}
%\endcases
$$
If $n$ is odd and $n+2k\equiv 1\mod4$ then
$$
\tilde b_{n,k}^J(\vec\alpha) = 
i^{\,-k} 2^{k-1} a_J^-(\vec\alpha),
\qquad
\tilde c_{n,k}^J(\vec\alpha) = 
i^k 2^{k-1} a_J^+(\vec\alpha).
$$
% $d_{n,k}^J(\alpha_1,\dots,\alpha_J)=$\newline
% $$%\qquad\qquad
% ={i\over4}(-2i^n)^{k+1} a_J^-(\alpha_1,\dots,\alpha_J)=
% \cases
% 	+i(-2i)^{k-1}\det A_J^- &\text{if $J\equiv 1\mod4$,}\\
% 	-i(-2i)^{k-1}\det A_J^+ &\text{if $J\equiv 3\mod4$.}
% \endcases
% $$
If $n$ is odd and $n+2k\equiv 3\mod4$ then
$$
\tilde b_{n,k}^J(\vec\alpha) = 
i^k 2^{k-1} a_J^+(\vec\alpha),
\qquad
\tilde c_{n,k}^J(\vec\alpha) = 
i^{\,-k} 2^{k-1} a_J^-(\vec\alpha).
$$
% $d_{n,k}^J(\alpha_1,\dots,\alpha_J) = $\newline
% $$%\qquad\qquad
% =-{i\over4}(-2i^n)^{k+1} a_J^+(\alpha_1,\dots,\alpha_J)=
% \cases
% 	+i(2i)^{k-1}\det A_J^+ &\text{if $J\equiv 1\mod4$,}\\
% 	-i(2i)^{k-1}\det A_J^- &\text{if $J\equiv 3\mod4$.}
% \endcases
% $$
\endproclaim

\demo{ Proof } By Lemma \lemCycPolA, it is sufficient to
compare the initial conditions which are computed in
Lemmas \lemDet, \lemDetBase, and \lemSignDetA. 
% Note also that if $n+2k\equiv 1\mod4$ then
% $$
% 	i^{nk}=i^{\,(2k+1)k}=i^{2k^2} i^k=(-1)^{k^2}i^k=(-1)^k i^k=(-i)^k.
% $$
% Analogously, if $n+2k\equiv-1\mod4$ then
% $i^{nk}=
% % =i^{(2k-1)k}=i^{2k^2} i^{-k} = (-1)^{k^2} i^{-k} = (-1)^k i^{-k}=(-i)^{-k}=
% i^k$.
\qed
\enddemo

%%%%%%%%%%%%%%%%%%%%%%%%%%%%%%%%%%%%%%%%%%%%%%%%%%%%%%%%%%%%%
%%%%%%%%%%%%%%%%%%%%%%%%%%%%%%%%%%%%%%%%%%%%%%%%%%%%%%%%%%%%%
%%%%%%%%%%%%%%%%%%%%%%%%%%%%%%%%%%%%%%%%%%%%%%%%%%%%%%%%%%%%%
%%%%%%%%%%%%%%%%%%%%%%%%%%%%%%%%%%%%%%%%%%%%%%%%%%%%%%%%%%%%%
%%%%%%%%%%%%%%%%%%%%%%%%%%%%%%%%%%%%%%%%%%%%%%%%%%%%%%%%%%%%%

\head \sectSign. Computation of the signatures
\endhead

\proclaim{ Lemma \lemSignDelta } Let $n$ and $k$ be positive integers.
The signature and the nullity of the link represented by the braid 
$\Delta^n\in B_{2k+1}$ are: 
$$
\split
\Sign\Delta^n &= 
\cases 	-nk(k+1) + (-1)^{(n-1)/2} 
		&\text{if $k\equiv n\equiv 1\mod2$},\\
	-nk(k+1) &\text{otherwise}.
\endcases
\\
\Null\Delta^n &= 
\cases	2k	&\text{if $n\equiv 0\mod4$},\\
	0	&\text{otherwise.}
\endcases
\endsplit
$$
\endproclaim

\demo{ Proof } Apply [\refNeumannSign]. \qed
\enddemo

\proclaim{ Proposition \propSignB } (a). Let $n$, $k$, $J$, and
$\alpha_1,\dots,\alpha_J$ be positive integers such
that $J\equiv n\mod2$. Let $b=b_{n,k}^J(\alpha_1,\dots,\alpha_J)$ 
be as in Section \sectDet.
When $n\equiv0\mod4$, we assume that $k=1$.
Then we have
{\rm(}$\eps_{n,k}$ is defined in {\rm(\defEpsnk))}
$$
\Null b = 
\cases
	1	&\text{if $J\equiv(2k-1)n\mod4$ and
			$\alpha_1=\dots=\alpha_J=1$},\\
	0	&\text{otherwise};
\endcases
$$
$$
\split
\Sign b - \Null b 
   \,=\, \Sign\Delta^n &+ (\alpha_1+\dots+\alpha_J) - J + 
   (-1)^k\cdot \Re i^{\,n-1}
\\
	= -nk(k+1) &+ (\alpha_1+\dots+\alpha_J) - J + \eps_{n,k} .
% \qquad\qquad\text{where}
\endsplit
$$
% $$
% \eps_{n,k}={1+(-1)^k\over2}\cdot\Re i^{\,n-1} =
% \cases
% 	(-1)^{(n-1)/2}	&\text{if $k+1\equiv n\equiv 1\mod2$,}\\
% 	0		&\text{otherwise}.
% \endcases
% $$

(b). Suppose that $n\equiv2\mod4$ and $k\ge1$.
Let $b=b_{n,k}^2(\alpha_0,0)$ for $\alpha_0\ge0$.
Then $\Sign b$ and $\Null b$ are computed by the formulas from the part (a)
with $J=0$ and the term ``\,$\alpha_1+\dots+\alpha_J$\!'' {\rm(}resp.
the condition ``\,$\alpha_1=\dots=\alpha_J=1$\!''\,{\rm)} replaced by
``\,$\alpha_0$\!'' {\rm(}resp. by ``\,$\alpha_0=1$\!''\,{\rm)}. 

(c). Suppose that $n\equiv0\mod4$ and $k = 1$.
Let $b=b_{n,k}^2(\alpha_0,0)$ for $\alpha_0\ge0$.
Then
$$
\Null b = 
\cases
	2	&\text{if $\alpha_1=\dots=\alpha_J=1$},\\
	1	&\text{otherwise};
\endcases
$$

\endproclaim

\proclaim{ Proposition \propSignC } (a). Let $n$, $k$, $J$, and
$\alpha_1,\dots,\alpha_J$ be positive integers such
that $J\equiv n\mod2$, $n\not\equiv 0\mod 4$, and $k\ge 2$.
Let $c=c_{n,k}^J(\alpha_1,\dots,\alpha_J)$ 
be as in Section \sectDet.
Then we have
{\rm(}$\eps'_{n,k}$ is defined in {\rm(\defEpsnkPrime))}
$$
\Null c = 
\cases
	1	&\text{if $J\equiv(2k+1)n\mod4$ and
			$\alpha_1=\dots=\alpha_J=1$},\\
	0	&\text{otherwise};
\endcases
$$
$$
\split
\Sign c - \Null c 
   \,=\, \Sign\Delta^n &+ (\alpha_1+\dots+\alpha_J) - J - 
   (-1)^k\cdot \Re i^{\,n-1}
\\
	= -nk(k+1) &+ (\alpha_1+\dots+\alpha_J) - J + \eps'_{n,k} .
% \qquad\qquad\text{where}
\endsplit
$$
% $$
% \eps'_{n,k}={1-3(-1)^k\over2}\cdot\Re i^{\,n-1} =
% \cases
% 	(-1)^{(n+1)/2}	&\text{if $k+1\equiv n\equiv 1\mod2$,}\\
% 	2(-1)^{(n-1)/2}	&\text{if $k\equiv n\equiv 1\mod2$,}\\
% 	0		&\text{otherwise}.
% \endcases
% $$
In particular, if $n$ is even and $k>1$ then $\Null c = \Null b$ and
$\Sign c = \Sign b$ for $b=b_{n,k}^J(\alpha_1,\dots,\alpha_J)$.

(b). Suppose that $n\equiv2\mod4$ and $k\ge2$.
Let $c=c_{n,k}^2(\alpha_0,0)$ for $\alpha_0\ge0$.
Then $\Null c = \Null b$ and
$\Sign c = \Sign b$ for $b=b_{n,k}^2(\alpha_0,0)$
%{\rm(}see Propsition {\rm\propSignC(b))}.

\endproclaim

\remark{ Remark } Computations show that if 
$n\equiv0\mod 4$ and $k\ge0$, then we have
$\Null b = 2(k-1) + \dots$, $\Null c = 2(k-1) + \dots$, 
$\Sign b - \Null b = -2(k-1)+\dots$, and
$\Sign c - \Null c = -2(k-1)+\dots$
where the dots stand for the corresponding expressions in 
Propositions {\propSignB} and \propSignC.
However, the method used in this paper is not sufficient
to prove this fact.
\endremark

\demo{ Proof of Proposition \propSignB } 
(The proof of Proposition {\propSignC} is analogous).
% The multiple induction using the values of $\Sign\Delta^n$
% and $\Null\Delta^n$ from Lemma {\lemSignDelta} for the base
% of induction and the skein relation (\skeinsigma) together
% with the signs of the determinants from Corollary {\corDnkJ}
% and Lemma {\lemSignDetA} for the induction step.
%
% We shall consider in details only the case when $n$ is odd.
% The case when $n$ is even, is similar. 

{\it Case 1}\; ($n$ {\it is odd\,}).
We shall use the induction
by $J$. Let us start with $J=1$.

We have 
$\tilde b^1_{n,k}(\alpha_1)=i^{-k}2^{k-1}a_1^-(\alpha_1)=
(-2i)^k(1+\alpha_1)$ if $n+2k\equiv1\mod4$ and 

\noindent
$\tilde b^1_{n,k}(\alpha_1)=i^k2^{k-1}a_1^+(\alpha_1)=
(2i)^k(1-\alpha_1)$
if $n+2k\equiv3\mod4$.
Hence, by Lemma \lemSkeinSigma, the signatures and the nullities are
as in Table 1 (we use that $b_{n,k}^1(0)\sim\Delta^n$).
Thus, the statement of the lemma holds for $J=1$.

\midinsert

\line {\hfill Table 1.  \hskip 8pt}
\smallskip
\moveleft 0pt\vbox {
\vbox{\tabskip=0pt 
%\offinterlineskip
\def \tablerule{\noalign {\hrule}}
\halign to 12.5cm {\strut#& \vrule# \tabskip=0.3em plus2em \hfil&
\hfil # \hfil & \vrule # &
\hfil # \hfil & \vrule # &
\hfil # \hfil & \vrule # &
\hfil #& \vrule#
\tabskip=5pt 

\cr\tablerule
&&      &&  $n+2k\equiv1\mod 4$ && $n+2k\equiv3\mod 4$
& \cr\tablerule
&& $\alpha_1$ && 
$\overset{\,}\to \;\;0\;\;\;1\;\;\;2\;\;\;3\;\;\;4\;\;\;5\;\;\;6\;\dots$ && 
		$\;\;0\;\;\;1\;\;\;2\;\;\;3\;\;\;4\;\;\;5\;\;\;6\;\dots$
& \cr\tablerule
&& ${\overset{\,}\to{\sign\big(\tilde b_{n,k}^1(\alpha_1)/
						\tilde b_{n,k}^1(0)\big)}}$ &&
		$\;  +\;\,  +\;    +\;\,  +\;    +\;\,  +\;    +  \dots$ &&
		$\;  +\;\;\,0\;\,  -\;\,  -\;    -\;\,  -\;    -  \dots$
& \cr\tablerule
&& ${\overset{\,}\to{\Sign b_{n,k}^1(\alpha_1)-\Sign\Delta^n}}$ &&
		$\;\;0\;\;\;1\;\;\;2\;\;\;3\;\;\;4\;\;\;5\;\;\;6\;\dots$ &&
		$\;\;0\;\;\;0\;\;\;0\;\;\;1\;\;\;2\;\;\;3\;\;\;4\;\dots$
& \cr\tablerule
&& ${\overset{\,}\to{\Null b_{n,k}^1(\alpha_1)}}$ &&
		$\;\;0\;\;\;0\;\;\;0\;\;\;0\;\;\;0\;\;\;0\;\;\;0\;\dots$ && 
		$\;\;0\;\;\;1\;\;\;0\;\;\;0\;\;\;0\;\;\;0\;\;\;0\;\dots$
& \cr\tablerule\cr}}
}
\endinsert

Now, 
suppose that we proved the statement of 
the lemma for smaller values of $J$.
Let us fix posivive integers $\alpha_2,\dots,\alpha_J$ and denote
$b(\alpha_1)=b_{n,k}^J(\alpha_1,\dots,\alpha_J)$, 
$\tilde b(\alpha_1)=\tilde b_{n,k}^J(\alpha_1,\dots,\alpha_J)$. 

We have 
$b(0) = b_{n,k}^{J-2}(\alpha_J+\alpha_2,\alpha_3,\dots,\alpha_{J-1})$.
By Corollary {\corDnkJ} and Lemma \lemSignDetA, the signs of the determinants 
are as in Table 2. Hence, by 
Lemma \lemSkeinSigma, the signature decrements and the nullities are
as in the next two lines of Table 2.

\midinsert

\line {\hfill Table 2.  \hskip 8pt}
\smallskip
\moveleft 0pt\vbox {
\vbox{\tabskip=0pt 
%\offinterlineskip
\def \tablerule{\noalign {\hrule}}
\halign to 12.5cm {\strut#& \vrule# \tabskip=0.1em plus2em \hfil&
\hfil # \hfil & \vrule # &
\hfil # \hfil & \vrule # &
\hfil # \hfil & \vrule # &
\hfil #& \vrule#
\tabskip=5pt 

\cr\tablerule
&&      && $n+2k\not\equiv J\mod 4$ &&   
& \cr
&&      && and $\alpha_j=1$ for $j>1$  && otherwise 
& \cr\tablerule
&& $\alpha_1$ && 
$\overset{\,}\to \;\;0\;\;\;1  \;\;\;2\;\;\;3\;\;\;4\;\;\;5\;\;\;6\;\dots$ && 
		$\;\;0\;\;\;1  \;\;\;2\;\;\;3\;\;\;4\;\;\;5\;\;\;6\;\dots$
& \cr\tablerule
&& ${\overset{\,}\to{\sign\big(\tilde b(\alpha_1)/\tilde b(0)\big)}}$ &&
		$\;  +\;\;\,0  \;\,  -\;\,  -\;    -\;\,  -\;    -  \dots$ &&
		$\;  +\;\,  -  \;    -\;\,  -\;    -\;\,  -\;    -  \dots$
& \cr\tablerule
&& ${\overset{\,}\to{\;\Sign b(\alpha_1)-\Sign b(0)+2}}$ &&
		$\;\;2\;\;\;2  \;\;\;2\;\;\;3\;\;\;4\;\;\;5\;\;\;6\;\dots$ &&
		$\;\;2\;\;\;1  \;\;\;2\;\;\;3\;\;\;4\;\;\;5\;\;\;6\;\dots$
& \cr\tablerule
&& ${\overset{\,}\to{\Null b(\alpha_1)}}$ &&
		$\;\;0\;\;\;1  \;\;\;0\;\;\;0\;\;\;0\;\;\;0\;\;\;0\;\dots$ && 
		$\;\;0\;\;\;0  \;\;\;0\;\;\;0\;\;\;0\;\;\;0\;\;\;0\;\dots$
& \cr\tablerule\cr}}
}
\endinsert

It remains to note that when $n$ is odd, 
$J\not\equiv n+2k$ iff $J\equiv (2k-1)n\mod 4$.

\smallskip
{\it Case 2}\; ($n\equiv 2\mod 4$). Similar to Case 1.

\smallskip
{\it Case 3}\; ($n\equiv 0\mod 4$).
For $\alpha_2\ne0$,
we have $\det b_{n,k}^J(0,\alpha_2)=0$ and $\det b_{n,k}(1,\alpha_2)\ne0$,
hence, $\Null b_{n,k}^J(0,\alpha_2)=1$ and 
$\Sign b_{n,k}^J(0,\alpha_2)=\Sign b_{n,k}^J(1,\alpha_2)$ by 
Lemma \lemSkeinSigma.
Since $\Null b_{n,k}^J(0,0)=2$, $\Null b_{n,k}^J(0,1)=1$, and
$\Null b_{n,k}^J(1,1)=0$, we have 
$\Sign b_{n,k}^J(1,1)=\Sign b_{n,k}^J(0,1)=\Sign b_{n,k}^J(0,0)=\Sign\Delta^n$.
The rest of the proof is similar to Case 1.
\qed
\enddemo

%%%%%%%%%%%%%%%%%%%%%%%%%%%%%%%%%%%%%%%%%%%%%%%%%%%%%%%%%%%%%
%%%%%%%%%%%%%%%%%%%%%%%%%%%%%%%%%%%%%%%%%%%%%%%%%%%%%%%%%%%%%
%%%%%%%%%%%%%%%%%%%%%%%%%%%%%%%%%%%%%%%%%%%%%%%%%%%%%%%%%%%%%
%%%%%%%%%%%%%%%%%%%%%%%%%%%%%%%%%%%%%%%%%%%%%%%%%%%%%%%%%%%%%
%%%%%%%%%%%%%%%%%%%%%%%%%%%%%%%%%%%%%%%%%%%%%%%%%%%%%%%%%%%%%

\head \sectProofMainTh. Proof of Theorem \MainTheorem
\endhead
We shall follow the scheme of the proof proposed in [\refOrevkovTop].
Let $A$ be a curve as in Theorem \MainTheorem. We associate to it 
a braid $b$ 
(the construction from [\refOrevkovTop; Sections 3.4 -- 3.5]
with $\Delta^n$ instead of $\Delta$). 
Then Murasugi-Tristram inequality imply
(see [\refOrevkovTop] for details)
$$
	\Null b + 1\ge |\Sign b\,| + m - e(b)
							\eqno(\ineqMT)
$$
where $e(b)$ is the {\it exponent sum} of $b$, i.e.
$e(b)=\sum k_j$ for $b=\prod \sigma_{i_j}^{k_j}$.

% Recall that $A_0$ is the union of the odd branch $\calJ$ and
% the double branches $O_1,\dots,O_{k-1}$ and the connected components
% of $\R A\setminus A_0$ are called the ovals. 
%
For an oval $v$, let $D_v$
be the component of $\pi_n^{-1}(\pi_n(v))\setminus(\calJ\cup E_n)$
which contains $v$.
Let $A'_\odd$ (resp. $A'_\even$) 
be the real (non-algebraic) smooth curve on $\R\Sigma_n$
which is obtained from $\R A$ by moving each odd (resp. even)
oval $v$ 
into the neighbouring component of $D_v\setminus A_0$.
Let $A''_\odd$ (resp. $A''_\even$) be the curve obtained from
$A'_\odd$ (resp. from $A'_\even$) by moving each oval $v$ into
the component of $D_v\setminus A_0$ which is the nearest to $\calJ$ under
the condition that the parity of $v$ is not changed
(see Figure {\figAodd}). Let us denote the corresponding braids by
$b'_\odd$, $b'_\even$, $b''_\odd$, and $b''_\even$ respectively.

\midinsert 
\epsfxsize 125mm
\centerline{\epsfbox{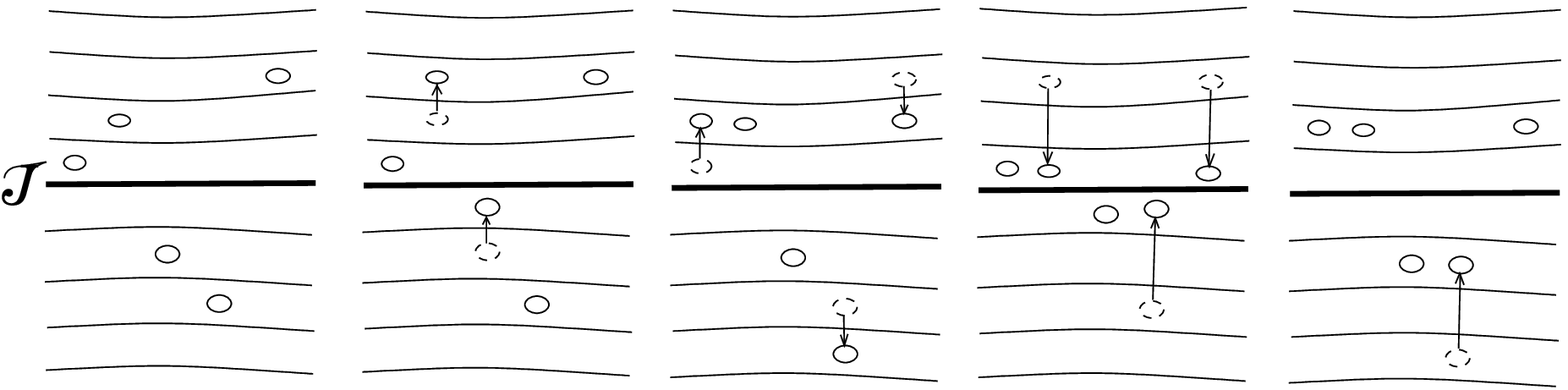}}
\centerline{\hbox to 3mm{}$\R A$ 
		\hbox to 16mm{} $A'_\odd$  \hbox to 14mm{} $A'_\even$
		\hbox to 14mm{} $A''_\odd$ \hbox to 14mm{} $A''_\even$
}
\botcaption{
             Fig.~\figAodd
          } 
\endcaption
\endinsert

Then we have $b''_\odd = b_{n,k}^J(\vec\alpha)$ and
$b''_\even = c_{n,k}^J(\vec\alpha)$ for some 
$\vec\alpha = (\alpha_1,\dots,\alpha_J)$ where $\alpha_j$ is
the number of ovals between two successive jumps over $\calJ$.
In particular, we have
$\alpha_1+\dots+\alpha_J=\lambda$. 
Thus, the signatures of $b''_\odd$ and $b''_\even$ can be
computed by Propositions {\propSignB} and {\propSignC}.
The condition $\lambda>J>0$ ensures that
$$
	\Null b''_\odd = \Null b''_\even = 0.		\eqno(\eqNullZero)
$$ 

When we pass from $b$ to $b'_\odd$ (resp. to $b'_\even$), we 
$\ODD$ (resp. $\EVEN$) times performe the replacement
$b_1\sigma_j^{-1} b_2 \longrightarrow b_1\tau \sigma_{j\pm1}^{-1}
\tau^{-1} b_2$ where 
$\tau = \tau_{j,j\pm1} = \sigma_{j\pm1}^{-1}\sigma_j$.
This is a composition of two band attachments:
$b_1\sigma_j^{-1} b_2 \longrightarrow b_1 b_2$ and
$b_1 b_2 \longrightarrow b_1\tau \sigma_{j\pm1}^{-1}\tau^{-1} b_2$
(it is clear that an inserting or a removing of $\sigma_j^{\pm1}$
is a band attachment). Hence, by Lemma \lemSkeinSigma(b), we have
$$
\matrix
	|\Sign b'_\odd - \Sign b\,| + |\Null b'_\odd - \Null b\,| 
						\le 2\ODD,
	\\
	\overset{\,}\to
	|\Sign b'_\even- \Sign b\,| + |\Null b'_\even- \Null b\,| 
						\le 2\EVEN.
\endmatrix
							\eqno(\ineqAodd)
$$

\midinsert 
\epsfxsize 100mm
\centerline{\epsfbox{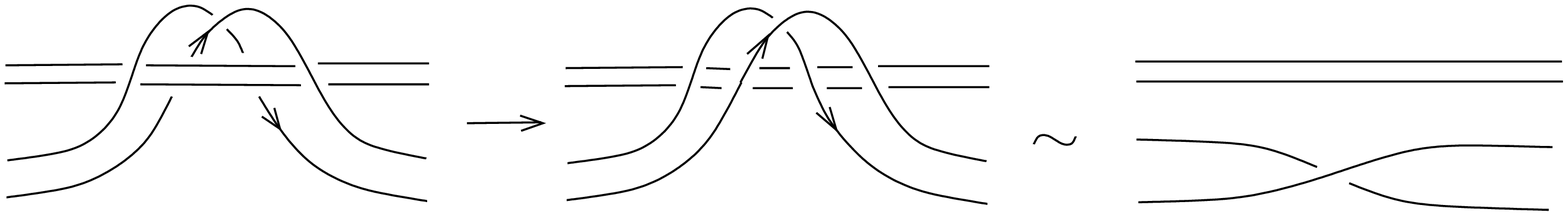}}
\botcaption{
             Fig.~\figModif. Passing from $b'_\odd$ to $b''_\odd$
          } 
\endcaption
\endinsert

When we pass from $b''_\odd$ to $b'_\odd$ (resp. 
from $b''_\even$ to $b'_\even$), we perform several times the
operation described in Lemma {\MainLemma}
(see Figure \figModif). Hence, by (\eqNullZero),
we have 
$$
\matrix
	\Sign b'_\odd =\Sign b''_\odd, & & \Null b'_\odd =0,
\\
	\overset{\,}\to
	\Sign b'_\even=\Sign b''_\even,& & \Null b'_\even=0.
\endmatrix
							\eqno(\eqAodd)
$$

Theorem {\MainTheorem} is a direct comination of
(\ineqMT), (\ineqAodd), (\eqAodd),
Propositions {\propSignB} and {\propSignC}, and the fact that
$e(b) = -(\alpha_1+\dots+\alpha_J) + e(\Delta^n) =
-\lambda + nm(m-1)/2$.

%%%%%%%%%%%%%%%%%%%%%%%%%%%%%%%%%%%%%%%%%%%%%%%%%%%%%%%%%%%%%
%%%%%%%%%%%%%%%%%%%%%%%%%%%%%%%%%%%%%%%%%%%%%%%%%%%%%%%%%%%%%
%%%%%%%%%%%%%%%%%%%%%%%%%%%%%%%%%%%%%%%%%%%%%%%%%%%%%%%%%%%%%
%%%%%%%%%%%%%%%%%%%%%%%%%%%%%%%%%%%%%%%%%%%%%%%%%%%%%%%%%%%%%
%%%%%%%%%%%%%%%%%%%%%%%%%%%%%%%%%%%%%%%%%%%%%%%%%%%%%%%%%%%%%
%%							   %%
%%							   %%
%%		    A P P E N D I C E S			   %%
%%							   %%
%%							   %%
%%%%%%%%%%%%%%%%%%%%%%%%%%%%%%%%%%%%%%%%%%%%%%%%%%%%%%%%%%%%%
%%%%%%%%%%%%%%%%%%%%%%%%%%%%%%%%%%%%%%%%%%%%%%%%%%%%%%%%%%%%%
%%%%%%%%%%%%%%%%%%%%%%%%%%%%%%%%%%%%%%%%%%%%%%%%%%%%%%%%%%%%%
%%%%%%%%%%%%%%%%%%%%%%%%%%%%%%%%%%%%%%%%%%%%%%%%%%%%%%%%%%%%%
%%%%%%%%%%%%%%%%%%%%%%%%%%%%%%%%%%%%%%%%%%%%%%%%%%%%%%%%%%%%%

\head 	Appendix A. Homogeneous skein systems and explicite formulas
	for the coefficients of $a_J^\pm$
\endhead

Let us introduce the following notation:
$\bar J=\{1,\dots,J\}$;
$[J]=\{ \{1,2\}, \{2,3\},\dots$, $\{J-1,J\}, \{J,1\}\}$;
$[[J]]=\{ A\subset[J] \mid \alpha\cap\beta = \varnothing \text{ for }
\alpha,\beta\in A \}$;
for $A\in[[J]]$, set $|A|=\bigcup_{\alpha\in A}\alpha$.
It is clear that $\card|A|$ is always even 
%(for a finite set $S$ we donote $\#S=\card S$).
For $k\equiv J\mod 2$, let us denote 
$\{A\in[[J]] \mid \card|A|=J-k\}$ by $[[J]]_k$ and set
$$
	f_{J,k}( x_1,\dots, x_J ) = 
	\sum_{A\in[[J]]_k} x_A
	\qquad\text{where}\quad
	x_A = \!\prod_{j\in\bar J\setminus|A|} x_j\,.
$$
It is clear that $f_{J,k}$ is a homogeneous polynomial of degree $k$.

Let us denote 
$\calF_k = \{f_{J,k}\}_{J=2,4,6,\dots}$ for an even $k$ and
$\calF_k = \{f_{J,k}\}_{J=1,3,5,\dots}$ for an odd $k$.
So, we have
$$
\xalignat1
%	&\calF_0=\{\,1,0,0,\dots\},	\qquad
%
	&\calF_1=\{\,x_1,\; x_1+x_2+x_3,\; 
			x_1+x_2+x_3+x_4+x_5,\,\dots\},
\\
	&\calF_2=\{\,%0,\; 
		x_1x_2,\;\;
		x_1x_2+x_2x_3+x_3x_4+x_4x_1,
\\
	&\qquad\qquad
		(x_1x_2+x_2x_3+x_3x_4+x_4x_5+x_5x_6+x_6x_1)+
		(x_1x_4+x_2x_5+x_3x_6),\,\dots\},
\\
	&\calF_3=\{\,0,\; x_1x_2x_3,\;\;
		x_1x_2x_3+x_2x_3x_4+x_3x_4x_5+x_4x_5x_1+x_5x_1x_2,
\\	
	&\qquad\qquad
	(x_1x_2x_3+x_2x_3x_4+x_3x_4x_5+x_4x_5x_6+x_5x_6x_7+x_6x_7x_1+x_7x_1x_2)+
\\	
	&\qquad\qquad
	(x_1x_2x_5+x_2x_3x_6+x_3x_4x_7+x_4x_5x_1+x_5x_6x_2+x_6x_7x_3+x_7x_1x_4),
	\,\dots\},
\\
	&\qquad\vdots
\endxalignat
$$

On can easily check that each $\calF_k$, $k\ge1$, is a 
skein system of cyclic polynomials.

\proclaim{ Proposition A.1 } 
If $J$ is even then
$$
	a_J^+(x_1,\dots,x_J) = 
	%\sum_{\smallmatrix k\ge 2\\k\equiv 0\mod 2\endsmallmatrix}
	\sum_{k=2,4,\dots,J}
	(2i)^k f_{J,k}(x_1,\dots,x_J)
	\quad\text{and}\quad
	a_J^- = -4 - a_J^+.
$$
If $J$ is odd then
$$
	a_J^+(x_1,\dots,x_J) = 2 +
	%\sum_{\smallmatrix k\ge 1\\k\equiv 1\mod 2\endsmallmatrix}
	\sum_{k=1,3,\dots,J}
	i(2i)^k f_{J,k}(x_1,\dots,x_J)
	\quad\text{and}\quad
	a_J^- = 4 - a_J^+.
$$
\endproclaim

\demo{ Proof } Follows from the fact that if $\{f_J\}$
is a skein system of cyclic polynomials then
$f_J=c_0 + \sum_{k\equiv J(2)} c_k f_{J,k}$ where
$c_0=f_J(0,\dots,0)$ and $c_k$ is the coefficient of
$x_1x_2\dots x_k$ in the polynomial $f_k$.
\enddemo

%%%%%%%%%%%%%%%%%%%%%%%%%%%%%%%%%%%%%%%%%%%%%%%%%%%%%%%%%%%%%%%%%%%
%%%%%%%%%%%%%%%%%%%%%%%%%%%%%%%%%%%%%%%%%%%%%%%%%%%%%%%%%%%%%%%%%%%
%%%%%%%%%%%%%%%%%%%%%%%%%%%%%%%%%%%%%%%%%%%%%%%%%%%%%%%%%%%%%%%%%%%
%%%%%%%%%%%%%%%%%%%%%%%%%%%%%%%%%%%%%%%%%%%%%%%%%%%%%%%%%%%%%%%%%%%
%%%%%%%%%%%%%%%%%%%%%%%%%%%%%%%%%%%%%%%%%%%%%%%%%%%%%%%%%%%%%%%%%%%
%%%%%%%%%%%%%%%%%%%%%%%%%%%%%%%%%%%%%%%%%%%%%%%%%%%%%%%%%%%%%%%%%%%
%%%%%%%%%%%%%%%%%%%%%%%%%%%%%%%%%%%%%%%%%%%%%%%%%%%%%%%%%%%%%%%%%%%

\head Appendix B. Generalized skein relations
\endhead
This appendix was written (as a section of the main body of the paper)
at the moment when I knew already Theorem {\thEN}
but did not understand that it immediately implies Corollary \corEN.
The identity (\genskeindetD) was used in the first version of the
proof of Lemma {\lemDet} to treat the case of non-fiberable links.
Now I do not know any application of these results but I decided to
keep them just because they look nice.

In this appendix we work with links presented in the form of closed 
braids. 
A link $L\subset S^3$ is presented 
by a closed braid with $m$ strings ($m$-braid)
if $p|_L$ is a covering $m$ of degree $m$ where $p$ is
the projection $S^3\setminus\ell = S^1\times\R^2 \to S^1$ for
some unknotted circle $\ell\subset S^3$.
We always suppose
the orientation of $L$ to be induced by the projection $L\to S^1$.
We use the language of braids just to siplify the notation. 
Everything can be reformulated for arbitrary link diagrams.
We always assume that $B_k\subset B_m$ for $k<m$ identifying $\sigma_j$ of
$B_k$ with $\sigma_j$ of $B_m$.

Set $\delta_k=\sigma_1\sigma_2\dots\sigma_{k-1}$.
Then $\Delta_k^2 = \delta_k^k$.
The skein relation (\skein) can be reformulated as follows: 
$\Omega_b + (t-t^{-1})\Omega_{b\delta} - \Omega_{b\delta^2}$
for any $b\in B_m$ and $\delta=\delta_2=\sigma_1$

% \proclaim{ Conjecture \conjGenSkein}
% For any $k$ there exist $n=n(k)$ and Laurent polynomials
% $c_{k,1}(t),$ $\dots,c_{k,n}(t)\in\Z[t,t^{-1}]$ such that
% $\sum_{j=0}^n c_{k,j}(t)\Omega_{b\delta_k^j}(t)=0$ for
% any $b\in B_m$, $m\ge k$.
% \endproclaim

\proclaim{ Proposition \propGenSkein }
  For any braid $b\in B_m$, $m\ge3$, one has
$$
%   \sum_{k=0}^4c_k(t)\Omega_{b\delta_3^k}(t) = 0     
   \Omega_b(t) + c_1\Omega_{b\delta}(t) + 
      c_2\Omega_{b\delta^2}(t) + c_3\Omega_{b\delta^3}(t) + 
           \Omega_{b\delta^4}(t) = 0                      \eqno(\genskein)
$$
where $\delta=\delta_3$, $c_1=c_3=-t^2+1-t^{-2}$,
      $c_2=-(t-t^{-1})^2$. In particular, for $t=i$,
$$
    \det b + 3\det(b\delta) + 4\det(b\delta^2)
        + 3\det(b\delta^3) + \det(b\delta^4) = 0.    \eqno(\genskeindet)
$$
\endproclaim

\demo{ Proof }
Starting with a link diagram, one can construct a Seifert surface $X$
using a standard algorithm based on so-called Seifert circles.
Applying this algorithm to a link presented as a closed braid
$b=\sigma_{i_1}^{\pm1}\dots\sigma_{i_n}^{\pm1}$
with $m$ strings, one obtains $m$ parallel equally oriented disks
and $n$ once-twisted ribbons where the $j$-th ribbon
connects the $i_j$-th disk to $(i_j+1)$-th one
(see Fig.~\figSeifertSurface\ for 
$b=\sigma_2\sigma_1^{-1}\sigma_2\sigma_2\sigma_1\in B_3$).

\midinsert 
\epsfxsize 65mm
\centerline{\epsfbox{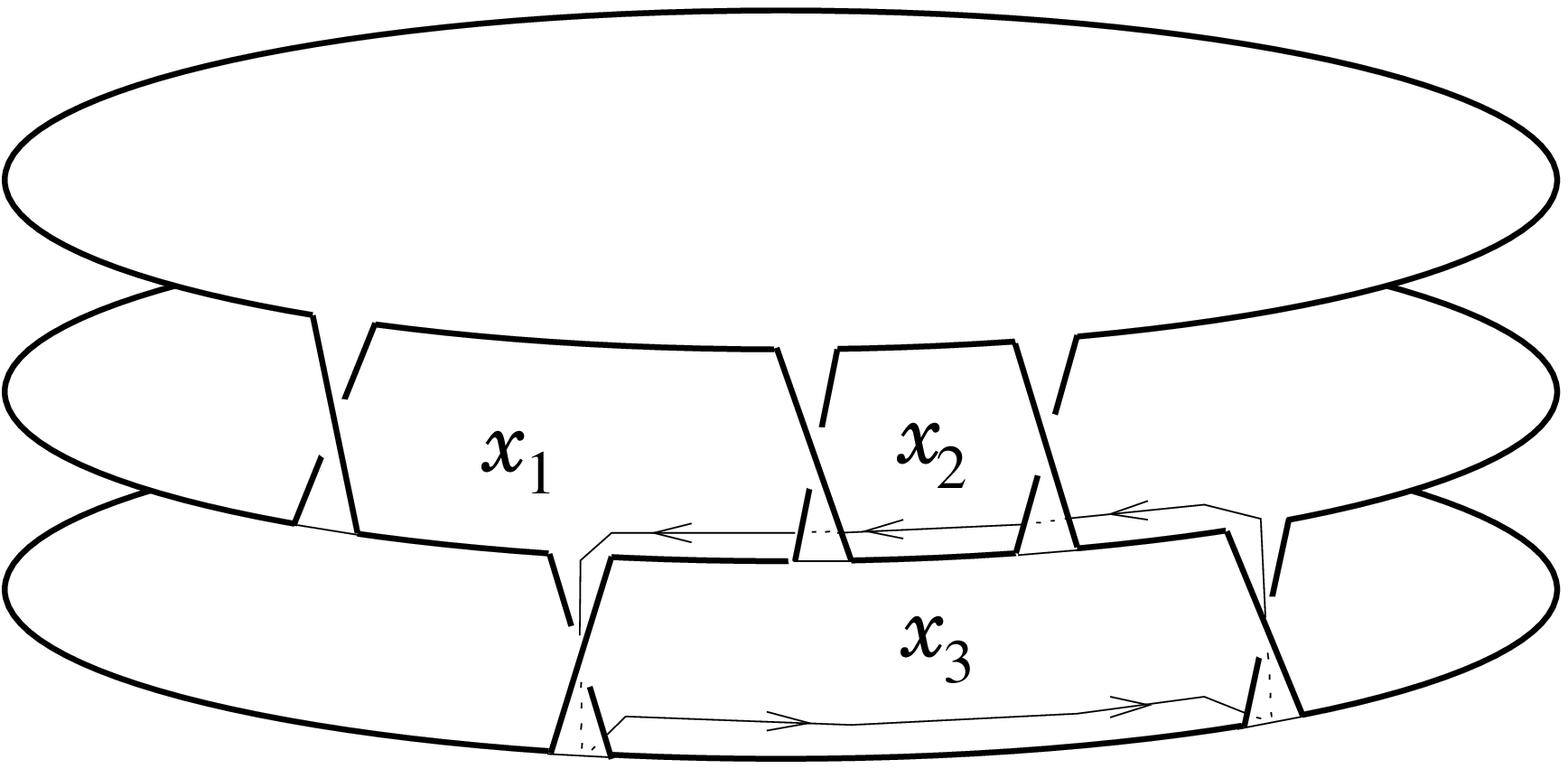}}
\botcaption{
             Fig.~\figSeifertSurface
          } 
\endcaption
\endinsert

As a base of $H_1(X)$ let us choose the cycles 
$x_1,\dots,x_{s}$, $s=n-m+1$ which
correspond to circuits in the positive direction around
the bounded regions of the projection of $b$ onto the plane
(cycles $x_1$, $x_2$, $x_3$ in Fig.~\figSeifertSurface).
Let $V=(v_{ij})$ be the corresponding Seifert matrix.
All the mutual positions of $x_\mu$ and $x_\nu$ providing 
$v_{\mu\nu}\ne0$ or $v_{\nu\mu}\ne0$ 
are shown schematically in Fig.~\figSeifertMatrix.

\midinsert
\epsfxsize 120mm
\epsfbox{ 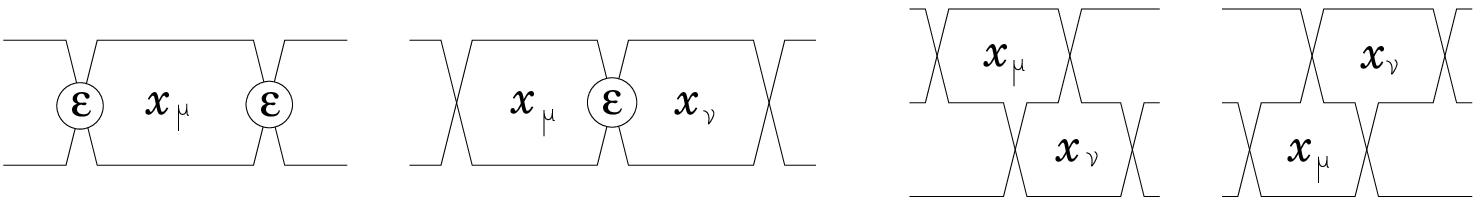 }
\smallskip
\vbox{\def\h{\hfil}\def\q{\qquad}\def\~{\tilde}
\halign{\q\quad#\h &\q\q\quad #\hfil &\q\q\quad #\hfil &\q\q #\hfil\cr
$v_{\mu\mu}=-\eps$&$v_{\mu\nu}=(\eps+1)/2$&$v_{\mu\nu}=0$&$v_{\mu\nu}=   -1$\cr
                  &$v_{\nu\mu}=(\eps-1)/2$&$v_{\nu\mu}=1$&$v_{\nu\mu}=\;\;0$\cr
}}
\botcaption{ Fig.~\figSeifertMatrix}
\endcaption
\endinsert

Aplying the above procedure to the braid $b_j:=b\delta^j$ 
we obtain the $(s+2j)\times(s+2j)$-matrix $V_j$ whose symmetrization
$V_j^t = t^{-1}V_j - tV^T$ has the form
$$
  V_j^t = 
   \left(\matrix
       V_0^t  &   U    &        &        & 0 \\
       U^*    &   W    &   B    &        &   \\
              &  B^*   &   A    & \ddots &   \\
              &        & \ddots & \ddots & B \\
	 0    &        &        &   B^*  & A 
   \endmatrix\right)
                                                 \eqno(\genskeinDemoA)
$$
$$
   \text{ where }\quad
   A=\left(\matrix
         t-t^{-1} &  -t^{-1}  \\
	   t      & t-t^{-1}
     \endmatrix\right),  \;\;
   B=\left(\matrix
        t^{-1}  &  0      \\
       -t       &  t^{-1}
     \endmatrix\right),   \;\;
   B^*=\left(\matrix
        -t &  t^{-1}  \\
	 0 & -t
     \endmatrix\right),
$$
$V_0^t$ is the $s\times s$ symmetrized Seifert matrix of $b=b_0$. 
Note, that $U$, $U^*$, and $W$ 
are $s\times 2$-, $2\times s$-, and $2\times 2$-matrices,
common for all $j$. 
Let us prove that 
$$
  \det V_0^t + c_1\det V_1^t + c_2\det V_2^t + c_3\det V_3^t + \det V_4^t = 0
                                                 \eqno(\genskeinDemoB)
$$
for matrices $V_j^t$ given by (\genskeinDemoA) where 
$V_0^t$, $U$, $U^*$, and $W$ 
are arbitrary fixed $s\times s$-, $s\times 2$-, $2\times s$-, 
and $2\times 2$-matrices
($c_j$ are the same as in (\genskein)).

Denote by $\widetilde V_j$ the $2j\times 2j$-matrix obtained from
the right lower $2j\times 2j$-minor of $V_j^t$ by replacing $W$
with a $2\times2$-matrix $\widetilde W$ with indeterminate 
entries $\widetilde w_{\mu\nu}$. Then 
$$
  \det\widetilde V_j = a_0^{(j)} + a_1^{(j)}\det\widetilde W + 
    a_{11}^{(j)}\widetilde w_{11} +
    a_{12}^{(j)}\widetilde w_{12} +
    a_{21}^{(j)}\widetilde w_{21} +
    a_{22}^{(j)}\widetilde w_{22}
$$
with $a_k^{(j)}, a_{\mu\nu}^{(j)}\in\Z[t,t^{-1}]$.
The staight forward computation shows that
$$
\xalignat3
a_0^{(2)}&=1&
a_0^{(3)}&=t^2-1+t^{-2}&
a_0^{(4)}&=t^4-t^2-t^{-2}+t^{-4}
\\
a_1^{(2)}&=t^2-1+t^{-2}&
a_1^{(3)}&=t^4-t^2-t^{-2}+t^{-4} &
a_1^{(4)}&=t^6-t^4+1-t^{-4}+t^{-6} 
\\
a_{12}^{(2)}&=t^3 &                    %b
a_{12}^{(3)}&=t^5-t^3 &
a_{12}^{(4)}&=t^7-t^5  
\\
a_{21}^{(2)}&=-t^{-3} &                    %a
a_{21}^{(3)}&=t^{-3}-t^{-5} &
a_{21}^{(4)}&=t^{-5}-t^{-7}  
\endxalignat
$$
and $a_{11}^{(j)}=a_{22}^{(j)}=t^{-2}a_{12}^{(j)}+t^2a_{21}^{(j)}$.
It is easy to see that $\det V_j^t$ has form
$$
  \det V_j^t = a_0^{(j)}\det V_0^t + a_1^{(j)}\det V_1^t + 
    a_{11}^{(j)}W_{22} +
    a_{12}^{(j)}W_{21} +
    a_{21}^{(j)}W_{12} +
    a_{22}^{(j)}W_{11}
                                                     \eqno(\genskeinDemoC)
$$
where $W_{\mu\nu}$ is the determinant of the matrix obtained from $V_1^t$
by deleting the row and the column contatining
the $(\mu,\nu)$-entry of $W$.

Substituting the expressions (\genskeinDemoC) for $\det V_2^t,\dots,\det V_4^t$
into the left hand side of (\genskeinDemoB), we obtain a linear
combination of $\det V_0^t$, $\det V_1^t$, and $W_{\mu\nu}$.
A straight forward computation shows that all the coefficients vanish.
\qed\enddemo

% \remark{ Remark } One can check that Conjecture \conjGenSkein\ also holds
% for $k=4$. In this case $n=8$ and
% $c_{4,0}=c_{4,8}=1$,
% $c_{4,1}=-c_{4,7}=t^3-t+t^{-1}-t^{-3}$,
% $c_{4,2}=c_{4,6}=-t^4+2t^2-2+2t^{-2}-t^{-4}$,
% $c_{4,3}=-c_{4,5}=t^5-2t^3+3t-3t^{-1}+2t^{-3}-t^{-5}$,
% $c_{4,4}=-t^6+t^4-3t^2+4-3t^{-2}+t^{-4}-t^{-6}$.
% \endremark

\proclaim{ Corollary \corGenSkein }
  For any braid $b\in B_m$, $m\ge3$, one has
$$
%   \sum_{k=0}^4c_k(t)\Omega_{b\delta_3^{2k}}(t) = 0     
   \Omega_b(t) + C_1\Omega_{b\Delta_3^2}(t) + 
      C_2\Omega_{b\Delta_3^4}(t) + C_3\Omega_{b\Delta_3^6}(t) + 
           \Omega_{b\Delta_3^8}(t) = 0                     \eqno(\genskeinD)
$$
where $C_1=C_3=-(t^3+t^{-3})^2$,
      $C_2=2(t^6+1+t^{-6})$. In particular, for $t=i$,
$$
    \det b - 2\det(b\Delta_3^4) + \det(b\Delta_3^8) = 0. \eqno(\genskeindetD)
$$
\endproclaim

\demo{ Proof }
Let (\genskein$_j$) be the result of substitution $b=\delta_3^j$
into (\genskein). Then (\genskeinD) is the sum of the identities
$(\genskein_0),\dots,(\genskein_8)$ multiplied by
$a_0,\dots,a_8$ respectively where
$a_0=a_8=1$, 
$a_1=a_7=t^2-1+t^{-2}$,
$a_2=a_6=t^4-t^2+1-t^{-2}+t^{-4}$,
$a_3=a_5=-t^4+t^2-2+t^{-2}-t^{-4}$,
$a_4=-2(t^2-1+t^{-2})$.
(we use here $\Delta_3^2=\delta_3^3$).
\qed
\enddemo

%%%%%%%%%%%%%%%%%%%%%%%%%%%%%%%%%%%%%%%%%%%%%%%%%%%%%%%%%%%%%%%
%%%%%%%%%%%%%%%%%%%%%%%%%%%%%%%%%%%%%%%%%%%%%%%%%%%%%%%%%%%%%%%
%%%%%%%%%%%%%%%%%%%%%%%%%%%%%%%%%%%%%%%%%%%%%%%%%%%%%%%%%%%%%%%
%%%%%%%%%%%%%%%%%%%%%%%%%%%%%%%%%%%%%%%%%%%%%%%%%%%%%%%%%%%%%%%
%%%%%%%%%%%%%%%%%%%%%%%%%%%%%%%%%%%%%%%%%%%%%%%%%%%%%%%%%%%%%%%
%%%%%%%%%%%%%%%%%%%%%%%%%%%%%%%%%%%%%%%%%%%%%%%%%%%%%%%%%%%%%%%
%%%%%%%%%%%%%%%%%%%%%%%%%%%%%%%%%%%%%%%%%%%%%%%%%%%%%%%%%%%%%%%

\enddocument